\documentclass[11pt]{article}
\usepackage{amsmath,amssymb,amsfonts}
\usepackage[dvips]{graphicx}

\newtheorem{theorem}{Theorem}
\newtheorem{lemma}{Lemma}		

\newcommand{\qed}{\hfill \ensuremath{\Box}}

\numberwithin{equation}{section}

\begin{document}

\title{\Large    Sufficient conditions for univalence of analytic functions }
\author{        D. Aharonov and U. Elias }
\date{}  
\maketitle


\begin{abstract}
The article discusses criteria for univalence of analytic functions in the unit disc. 
Various families of analytic functions depending on real parameters are considered. 
A unified method for creating new sets of conditions ensuring univalence is presented. 
Applying this method we are able to find several families of new sharp criteria for 
univalence.
\end{abstract}

\section{ Introduction }

The Schwarzian derivative  
    $ Sf = \left( {f''} / {f'} \right)' - \frac{1}{2} \left( {f''} / {f'} \right)^2 $ 
of an analytic locally univalent function plays an important role for finding sufficient 
conditions for univalence.  Z. Nehari  \cite{Nehari-1949} found conditions implying 
univalence expressed in terms of the Schwarzian derivative: 
If  $ |Sf| \le 2 ( 1-|z|^2)^{-2} $,  then  $ f $  is univalent in the unit disc  
$ \Delta = \{z , |z|<1 \} $.   Also if   $ |Sf| \le \pi^2 / 2 $,  the same conclusion 
follows.  For deriving his outstanding results Nehari used a useful connection between the 
zeros of solutions of linear second order differential equations and univalence
\cite{Nehari-1949}.  Later Pokornyi  \cite{Pokornyi}  stated without proof the condition  
$ |Sf| \le 4 ( 1-|z|^2)^{-1} $.  Nehari then proved this condition  \cite{Nehari-1954}.
In addition Nehari extended these results and proved a more general theorem
\cite{Nehari-1954},  \cite{Nehari-1979}  concerning criteria for univalence.
In his theorem he also investigated the sharpness of his conditions.

This pioneering works of Nehari opened a new line of research in geometric function theory.
One very significant development is due to Ahlfors and Weill  \cite{Ahlfors-Weill}, 
see also  \cite{Schober}.  In this paper the authors considered not only univalence 
criteria but in addition they were intersted in quasiconformal extension to the 
full complex plane   {\cal C}.   In another direction connection with harmonic 
mappings and minimal surfaces was made by Chuaqui, Osgood  and others
\cite{Chuaqui-1995},   \cite{Chuaqui-1998},   \cite{Chuaqui-2003},      
\cite{Chuaqui-2004},   \cite{Chuaqui-2007A},    \cite{Chuaqui-2007B}, \cite{Chuaqui-2008}, 
\cite{Chuaqui-2009}.

In the present paper we confine ourselves to finding new sufficient conditions for univalence 
by extending some of previous ideas.  In addition we deepen the  discussion of sharpness. 
In this connection we mention our paper  \cite{Aharonov-Elias}   which deals with singular 
Sturm comparison theorem.

As an example we show that our results lead to an improvement of the evaluation 
of the radius of univalence of the  ``error function".


\section{ Nehari's univalence criteria and it's sharpness  }

Nehari's pioneering work appeared in   \cite{Nehari-1949}.  
This work opened a fundamental line of research.  His idea was to use a 
 connection between the number of zeros of solutions of second order linear
differential equations in a given domain in the complex plane 
and univalence of the quotient of two independent solution of this equation: 
If  $ u(z), v(z) $  are two linearly independent solutions of a linear, homogeneous 
second order differential equation in a domain  $D$  and every 
linear combination  $ c_1 u(z) + c_2 v(z) $  has at most one zero in  $D$,  
then their quotient  $ f(z) = \ v(z) / u(z) $  is univalent in  $D$.

The well known Schwarzian derivative operator appeared naturally in his line
of reasoning and Nehari made use of the properties of this operator
to arrive at his sufficient conditions for univalence.
We now recall some of the basic properties of the Schwarzian Derivative ($ S$ 
by our notation):
$$
 Sf = \left( \frac{f''}{f'} \right)' - \frac{1}{2} \left( \frac{f''}{f'} \right)^2 
   = \frac{f'''}{f'} - \frac{3}{2} \left( \frac{f''}{f'} \right)^2 
$$ 
Some other commonly used notations are  $ Sf = \{f,z \},  \  Sf = S_f(z) $.

One basic property of the  $ S $   is:  Given a M\"obius map $ T = (az+b)/(cz+d) $, 
$ ad-bc \neq 0 $,  we have  $ S (T)(z) = 0 $.
Another useful property is for a composition of two functions  $ g \circ f $:
$$
	S ( g \circ f ) (z ) = ( S ( g ) \circ f(z) ) f'(z)^2 + S( f )(z) .
$$
If the above  $ f $  is in particular a M\"obius map  $ T $  then by  $ S(T) = 0 $, 
$$
	S ( g \circ T )(z) = ( S ( g ) \circ T(z) ) T'(z)^2  .
$$

Suppose we are given the linear differential equation
\begin{equation}						\label{eq:diffequ}
	    u'' + p(z) u = 0 ,
\end{equation}
where  $ p(z) $  is an analytic function in the unit disc  $ \Delta $.   
In addition, we shall need frequently the real valued differential equation
\begin{equation}						\label{eq:real-diffequ} 
	    y'' + p(x) y = 0 ,
\end{equation}
where  $ p(x) $  is a continuous function in the interval  $ -1 < x < 1 $.

Let  $ u(z), v(z) $  be any two linearly independent solutions of  (\ref{eq:diffequ}).  
Then it is a useful fact that 
\begin{equation}						\label{eq:S=2p}
	 S( v/u ) (z) = 2 p(z)    . 
\end{equation}
Nehari combined the above properties to arrive at his remarkable result

\noindent
{\bf Theorem A} (Nehari, \cite{Nehari-1954}) 
{\sl 
Suppose that  \\
(i) \ $ p(x) $  is a positive and continuous even function for $ -1 < x < 1 $,  \\
(ii) \ $ p(x) (1 - x^2) $  is nonincreasing for   $ 0 < x < 1 $,	\\ 
(iii) \ the differential equation   
			$ y''(x) + p(x)y(x) = 0 $  
has a solution which does not vanish in  $ -1 < x < 1 $.	\\
Then any analytic function  $ f(z) $  in  $ \Delta $   satisfying
\begin{equation}						\label{eq:Sf<p} 
	|  Sf(z) |  \le  2p( |z| ) 
\end{equation}
is univalent in the unit disc  $ \Delta $.
}

In what follows we use the term  ``Nehari's function"  to denote a positive even 
function  $ p(x) $  such that  $ p(x) ( 1 - x^2 )^2 $   is nonincreasing for 
$ 0 < x < 1 $.  See  \cite{Steinmetz}.

As Nehari pointed out already in  \cite{Nehari-1949},  the functions 
\begin{equation}						\label{eq:1-x^2^2}
	p(x) = ( 1 - x^2 )^{-2}, 	\qquad  p(x) = \pi^2 /4 ,
\end{equation} 
and the corresponding solutions  
   $ y(x) = (1 - x^2)^{1/2}, \ y(x) = \cos( \pi x / 2 ) $ 
of the respective equations  (\ref{eq:real-diffequ})  have all the needed properties 
to conclude the sufficient conditions for univalence in  $ \Delta $.
Soon after that,  E. Hille  \cite{Hille}  made the remarkable observation that 
the condition   $  | Sf | \le  2( 1 - |z|^2 )^{-2}  $   is sharp. 
In fact he showed more than that. Namely, if 
\begin{equation}						\label{eq:Hille}
	p(z) = (1 + \gamma^2 ) ( 1 - z^2 )^{-2} , 	\qquad  \gamma > 0, 
\end{equation}
then for two solutions  $  u(z), v(z) $  of the equation  $ u''(z) + p(z) u(z) = 0 $,
 the quotient   
	$ f(z) = v(z)/u(z) = \left( \dfrac{1+z}{1-z} \right)^{ i\gamma } $  
 satisfies
	$ |  Sf(z) |  \le 2(1 + \gamma^2 ) ( 1 - |z|^2 )^{-2}  $
but takes on the value 1 infinitely often in  $ \Delta $. 
On the other hand, if 
\begin{equation}						\label{eq:pi^2}
	p(z) = \frac{\pi^2} {4} + \delta , \qquad  \delta  > 0 ,
\end{equation}
the quotients of solutions of the corresponding differential equation  $ u'' + p(z) u = 0 $  
are finitely valent in   $ \Delta $.   Thus there is a clear difference between 
the two cases  (\ref{eq:Hille})  and  (\ref{eq:pi^2}).
In what follows we will pay attention to the question of sharpness and in 
particular we will investigate when the counterexamples are finitely or 
infinitely valent.

\noindent
{\bf Definition.} 
{\it
\  We shall say that analytic function  $ p(z) $  is
 ``self majorant"  in the unit disc  $ \Delta $  if 
\begin{equation}						\label{eq:self}
	     | p(z) | \le  p( |z| ) ,  \qquad  z \in \Delta .
\end{equation}
}

If all coefficients in the expansion of  $ p(z) $  around zero are positive,  
$ p(z) $  is automatically self majorant but the opposite is not true in general, 
as it follows by elementary examples.

The three first sufficient conditions implying univalence, namely
$$ 
 |Sf(z)| \le 2 ( 1-|z|^2)^{-2}, \qquad  |Sf(z)| \le \pi^2 / 2, \qquad  |Sf(z)| \le 4 ( 1-|z|^2)^{-1}
$$
correspond, to differential equations of type  (\ref{eq:real-diffequ})  where the 
respectively involved functions   
    $ p(z) = ( 1 - z^2)^{-2}, \ p(z) = \pi^2 /4,  \ p(z) = 2( 1 - z^2)^{-1} $ 
are all self majorant since the coefficients in their power series expansion are 
all positive.  This is not the case for Nehari's more 
general conditions  \cite{Nehari-1979}  for univalence in  $ \Delta $,    
\begin{align}
  |Sf(z)| &  \le 2(1+\mu) ( 1 - \mu |z|^2 ) ( 1 - |z|^2 )^{-2}, 
  \qquad\qquad		0 \le \mu \le 1, 		\label{eq:Nehari-mu-1} \\
  |Sf(z)| &  \le 2(1-\mu^2) ( 1 - |z|^2 )^{-2} + 2\mu( 2 + \mu ) ( 1 + |z|^2 )^{-2}, 
  \qquad	0 \le \mu \le 1,			\label{eq:Nehari-mu-2}
\end{align}
which correspond to the solutions $ y(x) = (1 - x^2)^{ (\mu+1)/2 } $  and 
$ y(x) = (1 - x^2)^{ (\mu+1)/2 } (1 + x^2)^{ -\mu/2 } $  of the respective equations
(\ref{eq:real-diffequ}).

In case the analytic function  $ p(z) $   is self majorant,  we have an 
additional significant information about the analytic differential equation
(\ref{eq:diffequ}): 
 If $ p(x) $ satisfies the assumptions of Theorem A and, in addition, 
$ p(z) $  is self majorant, then for two linearly independent solutions
$ u(z), v(z) $ of solutions of  (\ref{eq:diffequ}), the function  $ f(z) = v(z) / u(z) $ 
satisfies 
$$
		|Sf(z)| = | p(z) |  \le  p(|z|) , 
$$
hence  $ f(z) = v(z) / u(z)$  is itself univalent.  If, for a given solution  $ u(z) $  
of (\ref{eq:diffequ})  we choose the independent solution 
$$
	v(z) = u(z) \int^z \frac {dt} { u^2(t) } , 
$$	
then the quotient 
\begin{equation}					\label{eq:f=int-u^-2}
	     f(z)= \frac{ v(z) } { u(z) }  =  \int^z \frac {dt} { u^2(t) } 
\end{equation}
is itself univalent.  Hence in such a case univalent functions are naturally 
generated by the solutions of the differential equation  (\ref{eq:diffequ}).

In the spirit of N. Steinmetz  \cite{Steinmetz}, we define:

\noindent
{\bf Definition.} 
{\it
\  We shall say that the univalence criterion  (\ref{eq:Sf<p}) 
is sharp if  $  Sf(x)   \ge  2p(x) $  for  $ -1 < x < 1 $  and  
$ Sf(z)  \not\equiv  2p(z) $   imply that  $ f(z) $  is not univalent 
in  $ \Delta $.
}

Steinmetz proved the following result:

\noindent
{\bf Theorem B} (N. Steinmetz,  \cite[Section 3]{Steinmetz}) 
{\it 
If   $ p(z) $   is self majorant and  
	$ \int^1  y^{-2}(x) \, dx  = \infty $ 
for a positive solution  $ y(x) $  of  (\ref{eq:real-diffequ}), 
then the univalence criterion  (\ref{eq:Sf<p}) is sharp.  
Conversely, if  
	$ \int^1  y^{-2}(x) \, dx  < \infty $,  
then the univalence criterion  (\ref{eq:Sf<p})  is never sharp.
}

In a recent work  \cite{Aharonov-Elias}  the authors of the present paper 
proved a singular form of the classical Sturm comparison theorem for equations
of the form  $ y'' + p(x) y = 0 $.  We have shown there that if 
$ p(x) $  is continuous in  $ (-1, 1) $  but not necessarily bounded at
the endpoints  $ x = -1, 1 $,  then 
\begin{equation}					\label{eq:principal}
	\int_{-1} \frac{ dt } { y^2(t) } = \infty,	\qquad
	 \int^1   \frac{ dt } { y^2(t) } = \infty
\end{equation}
ensure the validity of the Sturm comparison theorem.
The assumptions  (\ref{eq:principal})  are closely related to that of \cite{Steinmetz},
which has to be satisfied in order to ensure the sharpness of the univalence condition.

One may ask why should not the definition of sharpness assume the inequality
$ | Sf(z) |  \ge  2|p(z)| $?   The following counter example explains our choice:
The differential equation
$$
	u'' - \frac{3} { (1 - z^2)^2 } \, u = 0 
$$
has the solutions  $ u_1(z) = (1+z)^{3/2} (1-z)^{-1/2},  \ u_2(z) = y_1(-z) $,
and the quotient of its solutions 
$$
	f(z) = \frac{ u_1(z) - u_2(z) }{ 4 u_2(z) } = \frac{z}{ (1-z)^2 }
$$
is Koebe's univalent function.  However, let
$$
	\frac{1}{2} | Sf(z) | = \left| \frac{ -3 } { (1 - z^2)^2 } \right|
		\ge   \left| \frac {1 + \gamma^2} { (1 - z^2 )^{2} } \right|
$$
and the quotient of solutions of Hille's equation  (\ref{eq:Hille})  with
$ 1 < 1 + \gamma^2  \le  3 $  is infinitely many valent.

We close this section with an example where  $ p(z) $  is not self majorant and
the criterion  (\ref{eq:Sf<p})  is not sharp.  Consider the inequality
\begin{equation}					\label{eq:Nehari-mu-3}
  |Sf(z)|  \le  6 ( 1 + |z|^2 )^{-2}
\end{equation}
which corresponds to the differential equation 
		$  y'' + \dfrac{3} { (1 + x^2)^2 } \, y = 0  $. 
(This is a particular case of  (\ref{eq:Nehari-mu-2})  with   $ \mu = 1 $).
Analyzing more carefully this condition we see that the condition  
(\ref{eq:Nehari-mu-3})  is not sharp. 
Indeed,  $ Sf(z) $  is an analytic function in the unit disc and because of the 
maximum principle,  $ \max | Sf(z) | $  cannot be achieved at any internal point of
the unit disc.  Consequently, it follows from  (\ref{eq:Nehari-mu-3})  that  
$ | Sf(z) | \le 3/2 $.  But this condition is not sharp, as we have the better 
sufficient condition  $ |Sf(z)|  \le  \pi^2 /2 $.  Note that the function 
$ p(z) = 3 ( 1 + z^2 )^{-2} $  is not self majorant.

The same example demonstrates that the assumption  $ | Sf(z) |  \le  2p(|z|) $   
is not necessary for the univalence of  $ f(z) = v(z) / u(z) $ .  Indeed,  
$ p(x) = 3 (1 + x^2)^{-2} $  satisfies assumptions (i)--(iii) of Theorem A,  
and the differential equation
$$
	u''(z) + \frac{3} { ( 1 + z^2 )^2 } \, u = 0 
$$
has the solutions  $ v(z) = z (1+z^2)^{-1/2}, \  u(z) = (1-z^2) (1+z^2)^{-1/2} $.   
However, for the function  $  \ f(z) = \dfrac{v(z)}{u(z)} = \dfrac{z}{ 1 - z^2 }  \ $, 
assumption  (\ref{eq:Sf<p}),  i.e., 
$$
	 | Sf(z) | = \left| \frac{6} { (1 + z^2)^2 } \right|
		\le   \frac {6} { (1 + |z|^2 )^{2} }  = 2p(|z|) 
$$
is false for nonreal values of  $z$  in the unit disc and Theorem A 
is not applicable.  Nevertheless this  $ f(z) $  is univalent in the unit disc.


\section{ Conditions for univalence depending on parameters }

Nehari's more general conditions  (\ref{eq:Nehari-mu-1}), (\ref{eq:Nehari-mu-2})  
for univalence    
\begin{align*}
  |Sf(z)| &  \le 2(1+\mu) ( 1 - \mu |z|^2 ) ( 1 - |z|^2 )^{-2}, 
  \qquad\qquad		0 \le \mu \le 1, 		 \\
  |Sf(z)| &  \le 2(1-\mu^2) ( 1 - |z|^2 )^{-2} + 2\mu( 2 + \mu ) ( 1 + |z|^2 )^{-2}, 
  \qquad	0 \le \mu \le 1, 
\end{align*}
may be viewed as families of conditions depending on one parameter.
Another family of conditions depending on a parameter were considered by Beesack
\cite{Beesack}, 
\begin{equation}					\label{eq:Beesack}
  |Sf(z)|  \le  \frac{\pi^2}{2} (1 - \lambda) + 2\lambda( 1 - \lambda )
  \left( \frac{\pi}{2} \tan \frac{ \pi |z| }{2}  -  \frac{2|z|}{ 1 - |z|^2 } \right)^2 ,
  \qquad   0  \le  \lambda  \le 1,  
\end{equation}
in  $ \Delta $.  It is generated by  
	$ y(x) = (1 - x^2)^\lambda  \cos^{1 - \lambda } (\pi x/2 )  $.
Here the coefficient  $ p(z) $  in the corresponding differential equation is 
self majorant and  (\ref{eq:Beesack})  is sharp due to Theorem B.

We now explain our method of finding families of conditions for univalence.
Our main idea is to consider a family of differential equations depending on 
one or two parameters. 
We take an even, bi-parametric function  $ u(z, \lambda, \mu) $  and generate through  
\begin{equation}				\label{eq:u''/u}
		u'' / u = -p(z, \lambda, \mu ) 
\end{equation} 
a differential equation 
\begin{equation}				\label{eq:p-lambda-mu}
	u'' + p(z, \lambda, \mu ) u = 0 .
\end{equation}
The idea is now to choose a suitable range of the two parameters 
$ \lambda, \mu $  to ensure that the corresponding  $ p(z, \lambda, \mu ) $   
will be a Nehari function and thus get for this range a set of criteria for 
univalence depending on the two parameters.  In fact, the number of parameters 
involved is not limited but the discussion becomes technically more involved.

In some cases we are lucky to have  $ p(z) $  be self majorant as well.
If this is not the case, some further restrictions on the parameters will imply
the self majorance property.

Taking all these in consideration, we now present the two families 
\begin{align*}
  u & = u(x,\lambda, \mu) = (1 - x^2)^\lambda G(x, \mu),
		\qquad  G(x) =G(-x),  \  G(x) \neq 0,  \quad
		  x \in (-1,1) , 
		\\
  u & = u(x,\lambda) = \cos(\pi x /2) G(x, \lambda),
		\qquad  G(x) =G(-x),  \  G(x) \neq 0,  \quad 
		 x \in (-1,1) .	
\end{align*}
In addition we also put for both families the restriction  $ u(1) = u(-1) = 0 $.  
This is a natural restriction since otherwise it is impossible to
get sharp conditions for univalence  (see  \cite{Steinmetz}  for more details).

We now have by  (\ref{eq:u''/u})

\begin{lemma}		
					\label{lemma:1} 
Let  $ u(x) =(1 - x^2)^\lambda G(x) $.  Then 
\begin{equation}					\label{eq:p-1-x2*G}
\begin{aligned}
    p(x, \lambda) =  &   -u'' / u			
  		  =  4\lambda( 1 - \lambda ) ( 1 - x^2 )^{-2}	\\
  &		  + 2\lambda ( 1 - x^2 )^{-1}				 
		  + 4\lambda x ( 1 - x^2 )^{-1} \left( {G'}/{G} \right) 
		  -  \left( {G'}/{G} \right)^2 - \left( {G'}/{G} \right)' .
\end{aligned}	
\end{equation}
\end{lemma}
The proof is elementary and details are omitted.  Some special cases of  
Lemma \ref{lemma:1}   are

\noindent  (a)
\begin{equation*}
\begin{aligned}
 u(x)  = &   (1 - x^2)^a \exp( \lambda x^2 ), 
 \\
 p(x, a, \lambda)  = &   4a(1-a) x^2 ( 1 - x^2 )^{-2} + 2a ( 1 - x^2 )^{-1} 
	       + 8a \lambda x ( 1 - x^2 )^{-1} -4\lambda x^2 - 2\lambda 		\\
    = &    (2a + 2a x^2 -4a^2 x^2) ( 1 - x^2 )^{-2} + 8a\lambda x^2 ( 1 - x^2)^{-1} 
			- (4\lambda^2 x^2 + 2\lambda) .
\end{aligned}
\end{equation*}
Note that  $ u(x) = (1 - x^2)^a \exp( \lambda x^2 )  $   is the solution 
of the initial value problem  $  u'' + p(x, a, \lambda)u = 0 $,  $ u(0)=1, \ u'(0) = 0 $.

\noindent  (b)
\begin{equation*}
\begin{aligned}
   u(x)  =  &  (1 - x^2)^\lambda \cos^\mu(\pi x /2),	\\
 p(x, \lambda, \mu) =  &   4\lambda( 1 - \lambda ) x^2 ( 1 - x^2 )^{-2} 
		+ 2\lambda ( 1 - x^2 )^{-1} + \mu\pi^2 / 4		\\
	 & + \mu(1 - \mu) \pi^2 \tan^2(\pi x /2) / 4 
	    -2 \mu \lambda \pi x  \tan(\pi x /2) ( 1 - x^2 )^{-1} .
\end{aligned}
\end{equation*}

\noindent  (c) 
\begin{equation*}
\begin{aligned}
  u(x)  = &  (1 - x^2)^a (1 + x^2)^{-b},  
  \\
 p(x, a, b)  = &  [ (2a -4ab) + x^2(2a + 4ab -4a^2) ] ( 1 - x^2 )^{-2} 	\\
	&  + [ (4ab + 2b) + x^2 (4ab -2b -4b^2) ] (1 + x^2)^{-2} .
\end{aligned}
\end{equation*}

Similarly to the above 

\begin{lemma}		
					\label{lemma:2} 
Let  \  $ u(x,\lambda) = \cos(\pi x /2) G(x, \lambda) $.   Then 
\begin{equation}					\label{eq:p-cos*G}
\begin{aligned}
  p(x, \lambda) = -u''(x) / u(x)  =  &  \pi^2 / 4 + \pi \left( {G'}/{G} \right) \tan(\pi x /2)  
  -  \left( {G'}/{G} \right)^2 - \left( {G'}/{G} \right)' .
\end{aligned}	
\end{equation}
\end{lemma}
Again, the elementary calculation is omitted.  Special cases of  Lemma \ref{lemma:2} are

\noindent  (d)  
\begin{equation*}
\begin{aligned}
  u(x) & =  \cos(\pi x /2) \exp( \lambda x^2 ), 
   \\
 p(x, \lambda) & = 2\lambda \pi x \tan(\pi x /2) + \pi^2 / 4 - (2\lambda + 4\lambda^2 x^2) .
\end{aligned}
\end{equation*}

\noindent  (e) 
\begin{equation*}
\begin{aligned}
 u(x) & =  \cos(\pi x /2) \exp( -\lambda \cos(\pi x /2) ),		\\
 p(x,\lambda) & =  (\pi^2 / 4)[ 1 + 2\lambda (1 - v^2) v^{-1} - \lambda^2 (1 - v^2) - \lambda v] ,
 		\qquad	v = \cos(\pi x /2).
\end{aligned}
\end{equation*}


\section{ Three sufficient conditions for univalence arising from Lemma \ref{lemma:1}  }

We start with the first condition arising from (a),  generated by 
$ u  = (1 - x^2)^a \exp( \lambda x^2 ) $.


\begin{theorem}
(A)  \    Let 
\begin{equation}				\label{eq:p-example1}		
 p(x) = (2a + 2a x^2 -4a^2 x^2) ( 1 - x^2 )^{-2} + 8a\lambda x^2 ( 1 - x^2)^{-1} 
	- (4\lambda^2 x^2 + 2\lambda) 
\end{equation}   
and let   $ \lambda, a $   satisfy 
\begin{equation}				\label{eq:example1-condition}  
  \frac{1 - 2a}{4}   \le  \lambda  
  \le  \frac{1}{2}   \left[ (1+2a) - (1+6a)^{1/2} \right], 
  \quad   
  \frac{1}{2}  \le  a  \le  1 .
\end{equation} 
  Then if  $ f(z) $  is an analytic function in   $ \Delta $  satisfying 
\begin{equation}				\label{eq:example1-bound}
	| Sf(z) | \le 2 p( |z| ) , \qquad  z \in \Delta , 				
\end{equation} 
it follows that  $ f(z) $  is univalent in  $ \Delta $.		\\
(B)  \  Let  $ p(x) $  be as in  (\ref{eq:p-example1})  and let  $ \lambda, a $  
satisfy the more restrictive conditions 
\begin{equation}				\label{eq:example1-condition-2} 
  \max \left\{   \frac{1 - 2a}{4}, a - \frac{1}{2} (6a)^{1/2}  \right\} 
  \le  \lambda  
  \le  \frac{1}{2}   \left[ (1+2a) - (1+6a)^{1/2} \right], 
  \quad   \frac{1}{2} < a \le 1 .
\end{equation} 
Then  $ p(z) $  is self majorant and the function 
\begin{equation}				\label{eq:4.5}
   f(z) = \int_0^z \frac{dt} { u^2(t) } 
   = \int_0^z \frac{dt} { (1 - t^2)^{2a} \exp( 2\lambda t^2 ) }	   	
\end{equation}
is an odd univalent function in  $ \Delta $.		\\ 
(C) \ Moreover, if the more restrictive conditions  (\ref{eq:example1-condition-2})  
are satisfied,  we conclude that the condition  (\ref{eq:example1-bound})  
is sharp.  \\
\end{theorem}

\noindent
{\bf  Proof of Theorem 1. }  \ 
We first note that the case  $ \lambda = 0 $  is elementary (also, it follows 
by continuity). 
Recalling how the condition for univalence was introduced (i.e., Lemma 1)  
and using Theorem A of Nehari, it is clear that the only thing  left is to check 
 that the function  $ p(x)( 1 - x^2 )^2 $  is positive and decreases on  $ [0,1] $. 
Hence we start with the expression  (\ref{eq:p-example1})  for  $ p(x) $  and  
$ a, \lambda $  real.  
Our aim is to find conditions on the parameters  $ a, \lambda $  that will ensure 
positivity and monotonicity of  $ p(x)( 1 - x^2 )^2 $.

From  $ p(0) = 2a - 2\lambda \ge 0 $  it follows at once that we need 
\begin{equation}					\label{eq:a>lambda} 	
		a \ge \lambda .
\end{equation} 
Also  $ \lim_{ x \to 1^- }p(x)( 1 - x^2 )^2 \ge 0  $  implies the restriction 
\begin{equation}					\label{eq:a<1} 		
		4a(1 - a) \ge 0 ,
\end{equation}
i.e.,  $ 0 \le a  \le 1 $.

For  $ p(x) $  in  (\ref{eq:p-example1})  it will be convenient to denote 
$$ 
	\varphi(x) :=  p(x)( 1 - x^2 )^2 = A + B x^2 + C x^4 + D x^6 ,
$$
 and with  $ t = x^2 $,  let   $  \psi(t) = A + Bt + C t^2 + D t^3 $. 
Then by direct computation
\begin{equation}					\label{eq:ABCD} 	
\begin{aligned}
	& 	A = 2(a-\lambda), \quad
		B = 2(a - 2 a^2 + 4a\lambda - 2\lambda^2 + 2\lambda),	\\
	&	C = -2\lambda( 4a  - 4\lambda + 1),	\quad
		D = -4 \lambda^2 .  
\end  {aligned}
\end{equation}
In order to ensure that   $ p(x)( 1 - x^2 )^2 $  decreases, we shall require that 
$ \psi'(t) \le 0 , \ 0 \le t \le 1 $.   For this aim it will enough to prove  that 
\begin{equation}						\label{eq:psi'(t)}  	
	\psi'(0) \le 0,   \quad    \psi'(1) \le 0   
	\quad 	\hbox{\rm and }	\quad   \psi''(t) \neq 0 , \quad  0 \le t \le 1 . 
\end{equation}
From  $ \psi'(0) = B = 2(a - 2 a^2 + 4a\lambda - 2\lambda^2 + 2\lambda) \le 0 $ 
  we have
\begin{equation}					\label{eq:4.13} 	   	
	2 \lambda^2 - 2\lambda (2a+1) + (2a^2 -a)  \ge  0 .
\end{equation}
Denote   $ \lambda_{1,2} = [ (1+2a) \pm (1+6a)^{1/2} ] / 2 $.  
It follows from  (\ref{eq:a>lambda})  that  $ \lambda \ge \lambda_2 $   is impossible 
and we are left with the condition
\begin{equation}					\label{eq:example1-lambda1} 	
	\lambda  \le  \lambda_1 = [ (1+2a) - (1+6a)^{1/2} ] / 2 .
\end{equation}
Note that in view of  $ a \ge 0 $,   (\ref{eq:a>lambda})  is contained 
in (\ref{eq:example1-lambda1}),  since   $ [ (1+2a) - (1+6a)^{1/2} ] / 2 \le a $ 
indeed holds if $ a \ge 0 $.

Since   $ \psi'(1) = B + 2C + 3D \le 0 $   simplifies to  $ a( 1 - 2a -4\lambda) \le 0 $ 
and  $ a \ge 0 $,   we get
\begin{equation}					\label{eq:4.16} 	   
	\lambda \ge (1 - 2a) / 4 .
\end{equation}
So far we have the restrictions 
\begin{equation}					\label{eq:4.17} 	   
	 (1 - 2a) / 4 \le \lambda \le [ (1+2a) - (1+6a)^{1/2} ] / 2,
	\quad   0 \le a  \le 1 , 
\end{equation}
which imply   $ \psi'(0) \le  0,   \    \psi'(1)  \le  0  $.
The inequality  $ (1 - 2a)/4  \le  [ (1+2a) - (1+6a)^{1/2} ] / 2 $ 
in  (\ref{eq:4.17})  also implies that  $ a \ge \frac{1}{2} $.   Thus we
have in  (\ref{eq:example1-condition})  the more restrictive condition 
$ a \ge \frac{1}{2} $  rather then  $ a \ge 0 $.

Finally, to complete the proof of  (\ref{eq:psi'(t)}),  we demand that     
	$ \psi''(t) \neq 0 , \ 0 \le t \le 1 $. 
If this is not the case, there exists  $ t_0,  \ 0 \le t_0 \le 1 $,  
such that  $ \psi''(t_0) = 2C + 6D t_0 = 0 $,  and by   (\ref{eq:ABCD}), 
$  t_0 = -2C / 6D = -\lambda(4a  - 4\lambda + 1) / 6 \lambda^2 $.
Then  $ 0 \le t_0 \le 1 $  is equivalent to 
\begin{equation}					\label{eq:4.19} 	   
	\lambda (4\lambda -4a - 1) \ge 0,  \qquad   \lambda (2\lambda + 4a + 1) \ge 0 .
\end{equation}
First consider the case   $ \lambda > 0 $.  We then get  $ \lambda  \ge a + 1/4 $,  
which contradicts  $ a \ge  \lambda $   in  (\ref{eq:a>lambda}).
Assuming now   $ \lambda < 0 $.  This implies  $ 2\lambda + 4a + 1 \le 0 $   
which contradicts   $ \lambda \ge (1 - 2a) / 4 $   in  (\ref{eq:4.16}), 
as $ a \ge 0 $.  These contradictions verify that  $ \psi''(t) \neq 0 $
for  $ 0 \le t \le 1 $  and complete the proof part (A).

We now turn to the proof of part (B).  Our aim is to find conditions that imply 
the self majorance property. For this aim it is enough to ensure positivity 
of the coefficient of every  $ x^{2n} $  in the expansion of the even function $ p(x) $.
From  (\ref{eq:p-example1})  with  $ t = x^2 $, 
\begin{equation}					\label{eq:p(x=t^2)} 	   
	p(x)/2 = (a + at -2a^2 t) \sum_{n=0}^\infty (n+1) t^n 
		+ 4a\lambda \sum_{n=1}^\infty t^n - (2\lambda^2 t + \lambda) .
\end{equation}
We shall require that the coefficient of each  $ t^n $  in  (\ref{eq:p(x=t^2)}) 
is positive.

The coefficient of  $ t^0 $  is   $ a - \lambda \ge 0 $.
For the coefficient of  $ t^1 $  to be positive,  we demand that
\begin{equation}					\label{eq:4.t^1} 	   
	3a - 2a^2 + 4a\lambda -2 \lambda^2 \ge 0 ;
\end{equation}
and for the coefficient of  $ t^n, n \ge 2 $,   we demand that
\begin{equation}					\label{eq:4.t^n} 	   
	(n+1)a + n(a - 2a^2) + 4a\lambda \ge 0 .
\end{equation}
Observe that  (\ref{eq:4.t^1})  implies  (\ref{eq:4.t^n}).  Indeed, 
$$
	(n+1)a + n(a - 2a^2) + 4a\lambda  \ge 3a - 2a^2 + 4a\lambda - 2 \lambda^2 
$$
is equivalent to 
$ (2n-2)(a - a^2) \ge  - 2\lambda^2 $,  which is true for all  $ n \ge 2 $  and  
$ 0 \le a \le 1 $.   So we have to take care only of  (\ref{eq:4.t^1}),  i.e.,   
	$  \ 2 \lambda^2 - 4a\lambda + (2 a^2 -3a) \le 0  \ $.
Since  $ \lambda \le a $,  we are left with the condition  
\begin{equation}					\label{eq:4.23} 	   
	a \ge \lambda  \ge a - (6a)^{1/2} / 2 .
\end{equation}
Summing up we conclude that the conditions ensuring  $ | p(z) | \le p( |z| ) $ 
are exactly as claimed in  (\ref{eq:example1-condition-2}).

The relation  (\ref{eq:4.5})  is established in  (\ref{eq:f=int-u^-2}). 
The divergence of the integral  (\ref{eq:4.5})  at  $ z=1 $  for  $ a \ge 1/2 $ 
implies the sharpness claimed in part (C)  according to Theorem B.
Now the proof of Theorem 1 is complete.

We end up with two figures:  Figure 1  illustrates the conditions 
(\ref{eq:example1-condition}) of part (A)  and Figure 2 demonstrates the 
(more restrictive) condition (\ref{eq:example1-condition-2}) of part (B). 
For Figure 2  observe that    $ a - (6a)^{1/2} / 2 = (1 - 2a) / 4  \  $  
implies   $  \  a = 1/2 + \sqrt{2}/3 = 0.9714...  $.
\qed

\begin{center}
    \includegraphics[scale=0.8]{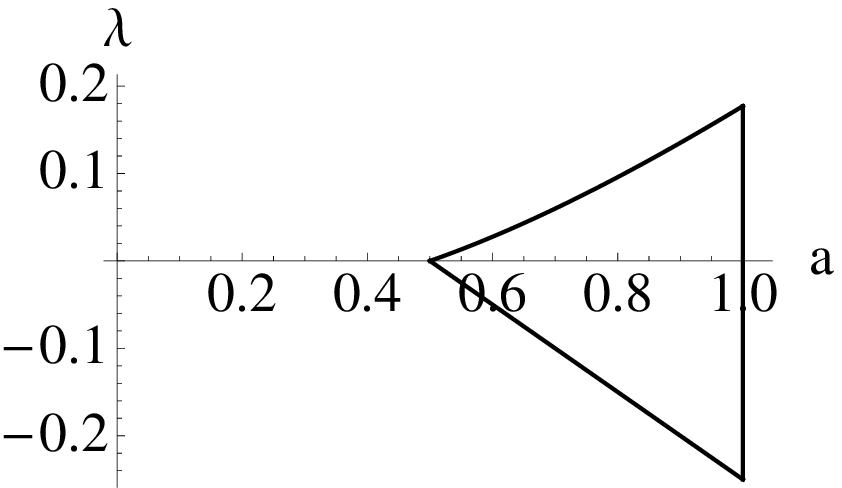}
\end{center}
\centerline{ Figure 1. The domain defined by  (\ref{eq:example1-condition})  }

\bigskip\bigskip
\begin{center}
    \includegraphics[scale=0.8]{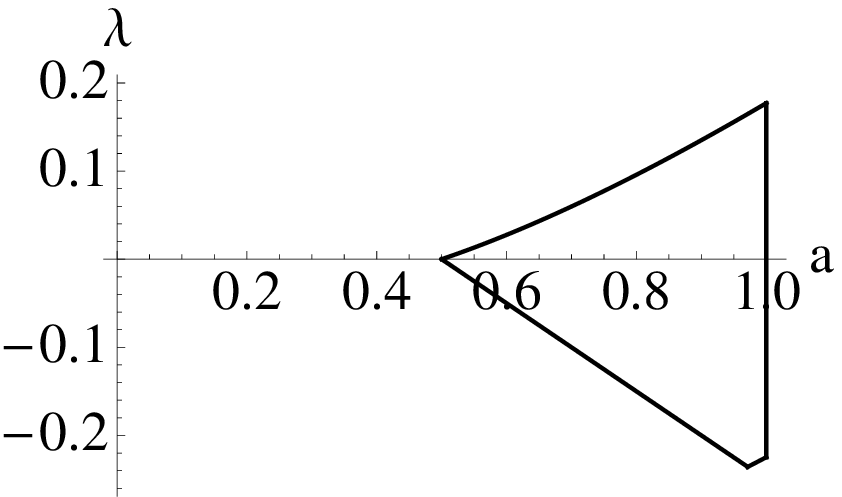}
\end{center}
\centerline{ Figure 2.  The domain defined by  (\ref{eq:example1-condition-2}) }

\bigskip


\noindent
{\bf Remark.}  
Recall that  $ u(z) = (1 - z^2)^a \exp( \lambda z^2 ) $  is an even function and  $ f(z) $
in  (\ref{eq:4.5})  is odd.  It is a well known classical fact that for every odd 
univalent function  $ f(z) = z + a_3 z^3 + \ldots $  in the unit disc, 
the function  $ g(z) $  defined by   $ \left( g(z^2) \right)^{1/2} = f(z) $  
is a normalized univalent function in the unit disc  \cite[p. 213]{Nehari-book}.
Hence, Theorem 1, (B)  may generate an additional univalent function.
This observation applies to the other four theorems as well.

The next example is generated by the function 
	$ u(x)  = (1 - x^2)^\lambda \cos^\mu(\pi x /2) $,  
which is introduced in case (b) of Lemma 1. 


\begin{theorem} 
Let
\begin{equation}				\label{eq:p-example2}	
\begin{aligned}
 p(x) = &  4\lambda( 1 - \lambda ) x^2 ( 1 - x^2 )^{-2} + 2\lambda ( 1 - x^2 )^{-1} 
        +\mu\pi^2 / 4		\\
	&  + \mu(1 - \mu) \pi^2 \tan^2(\pi x /2) / 4 
		-2 \mu \lambda \pi x  \tan(\pi x /2) ( 1 - x^2 )^{-1} 	
\end  {aligned}
\end{equation}
and let   $ \lambda, \mu $   satisfy 
\begin{equation}				\label{eq:example2-condition}  
  \lambda \ge 0,	\quad   \mu \ge 0, \quad   1/2 \le \lambda + \mu \le 1,
  \quad (1 - \mu)/2 \le \lambda \le 1 - \mu .
\end{equation} 
  Then if  $ f(z) $  is an analytic function in   $ \Delta $  satisfying 
\begin{equation}				\label{eq:4.26}		
	| Sf(z) | \le 2 p( |z| ) , \qquad  z \in \Delta , 	
\end{equation} 
it follows that  $ f(z) $  is univalent in  $ \Delta $.		\\
Also
\begin{equation}				\label{eq:4.27}		
   f(z) = \int_0^z \frac{dt} { u^2(t) } 
   = \int_0^z \frac{dt} { (1 - t^2)^{2\lambda} \cos^{2\mu}( \pi t/2 ) }	
\end{equation}
is an odd univalent function in  $ \Delta $.
Moreover the condition  (\ref{eq:4.26})  is sharp.
\end{theorem}

\noindent
{\bf Proof of Theorem 2.}  \ 
Using an identical argument as in Theorem 1 we proceed to ensure that the 
function   $ \varphi(x) = p(x) (1 - x^2)^2 $  is nonincreasing for  $ 0<x<1 $.
We start with the expression  (4.24) for  $ p(x) $   with   $ \lambda, \mu $  
real.  Our aim is to find the restrictions on  $ \lambda, \mu $  that will imply 
the monotonicity of   $ \varphi(x) $.

It will be convenient to denote
\begin{equation}				\label{eq:G=(1-x^2)tan}	 
	G(x) := (1 - x^2) \tan(\pi x /2) .
\end{equation}
With this notation, 
\begin{equation}				\label{eq:phi-example2}	
  \varphi(x) = 4\lambda ( 1 - \lambda ) x^2 + 2\lambda ( 1 - x^2 ) 
	+ \mu\pi^2 ( 1 - x^2 )^2 / 4 + \mu(1 - \mu) \pi^2 G^2(x) / 4 
	-2 \mu \lambda \pi x G(x) .	
\end{equation}

We start with some elementary considerations.  To ensure that  $ \varphi(x) $  
is positive and nonincreasing for  $ 0 < x < 1 $,  we must have in particular that 
$ \varphi(1) \ge 0 $  and  $ \varphi'(1) \le 0 $.  Thus we require that 
\begin{equation}				\label{eq:phi(1)-example2}	
	\varphi \Big|_{x=1^-} = 4( \lambda + \mu ) (1 -\lambda - \mu) \ge 0 ,
\end{equation}
i.e.,  $ 0 \le \lambda + \mu \le 1 $.
If  $ \lambda + \mu = 0 $, then from   $ u  = (1 - x^2)^\lambda \cos^\mu(\pi x /2) $ 
it easily follows that  $ u(1)=u(-1) $  are different from zero, which contradicts 
our restrictions on  $u$.  Hence let  $ 0 < \lambda + \mu \le 1 $.

We have from  (\ref{eq:phi-example2}) that 
\begin{equation}				\label{eq:phi'}	
	\varphi'(x) = 8\lambda ( 1 - \lambda ) x - 4\lambda x 
	-  \mu\pi^2 x ( 1 - x^2 ) + \mu(1 - \mu) (\pi^2/ 4 ) 2G(x) G'(x)  
	-2 \mu \lambda \pi (x G'(x) + G(x) ) .
\end{equation}
By direct calculation, 
\begin{equation}				\label{eq:G(1)}		
	G \Big|_{x=1^-} = 4 / \pi,	\qquad	G' \Big|_{x=1^-} = 2 / \pi ,
\end{equation}
and after some more calculations, we get  
	$ \varphi'(1) = 4( \lambda + \mu ) (1 -2\lambda - \mu) $. 
Hence we require 
\begin{equation}				\label{eq:phi'(1)-example2}	
	\varphi'(1) = 4( \lambda + \mu ) (1 -2\lambda - \mu) \le 0 .
\end{equation}
Summing up, we get from  (\ref{eq:phi(1)-example2})  and (\ref{eq:phi'(1)-example2}) 
the conditions 
\begin{equation}				\label{eq:4.34}		
	1 - \lambda - \mu \ge 0,    \quad	1 -2\lambda - \mu \le 0. 
\end{equation} 
Consequently  $ \lambda \ge 0 $  and   $ \mu \le 1 $  are necessary for  $ \varphi(x) $ 
to be nonincreasing.

From now and on we assume that  $ \mu > 0 $.  The case  $ \mu \le 0 $  will be 
mentioned later, in a remark at the end of the proof.

We proceed to show that under the restrictions above,  $ \varphi(x) = p(x)(1 - x^2)^2 $  
is indeed positive and nonincreasing.  For this purpose we rewrite  $ p(x) $  of  
(\ref{eq:p-example2})  as 
\begin{equation*} 
\begin{aligned}
 p(x) = &  4\lambda( 1 - \lambda ) x^2 ( 1 - x^2 )^{-2} + 2\lambda ( 1 - x^2 )^{-1} 
          + \mu\pi^2 / 4		\\
&  + \mu(1 - \mu) 
 	\left( 
		(\pi/2)\tan(\pi x /2) - \frac{2 \lambda}{1-\mu} \frac{x} { 1 - x^2 }  
 	\right)^2
	-\frac{4 \mu \lambda^2}{1-\mu} x^2 ( 1 - x^2 )^{-2} 	
\end  {aligned}
\end{equation*}
and denote 
\begin{equation}				\label{eq:L(x)}		
    L(x) = (\pi/2) \tan ( \pi x / 2 ) - \frac{2\lambda }{1-\mu} \, \frac{x}{1 - x^2 } .
\end{equation} 
With this notation,
\begin{equation*}
 p(x) =   \frac{ 4\lambda( 1 - \lambda -\mu )} { 1-\mu } x^2 ( 1 - x^2 )^{-2} 
        + 2\lambda ( 1 - x^2 )^{-1} + \mu\pi^2 / 4 
	+ \mu(1 - \mu) L^2(x)  ,
\end{equation*}
\begin{equation*}
\varphi(x) = \frac{ 4\lambda( 1 - \lambda -\mu )} { 1-\mu } x^2  
        + 2\lambda ( 1 - x^2 )
        + ( \mu\pi^2 / 4 ) ( 1 - x^2 )^2 	\\
	+ \mu(1 - \mu) \left[ ( 1 - x^2 ) L(x) \right]^2 
\end{equation*}
and
\begin{equation*}
\varphi'(x) = \frac{ 2\lambda( 1 - 2\lambda -\mu )} { 1-\mu } 2x 
        -  \mu\pi^2 x( 1 - x^2 ) 	\\
	+ 2\mu(1 - \mu) \left[ ( 1 - x^2 ) L(x) \right] \left[ ( 1 - x^2 ) L(x) \right]' .
\end{equation*}
As  $ 1 - \lambda -\mu \ge 0 $  by  (\ref{eq:4.34}),  $ \lambda \ge 0 $  and 
$ 0 \le \mu \le 1 $,  it is clear that  $ \varphi(x) \ge 0 $. 
But we also have  $ 1 - 2\lambda -\mu  \le 0 $  by  (\ref{eq:4.34}).  So, in order 
to establish  $ \varphi'(x) \le 0 $,  it is left to show only that
\begin{equation}				\label{eq:L'-ineq}			
  2\mu(1 - \mu) \left[ ( 1 - x^2 ) L(x) \right] \left[ ( 1 - x^2 ) L(x) \right]' 
	\le	\mu\pi^2 x( 1 - x^2 ) 
\end{equation}
We postpone the proof of (\ref{eq:L'-ineq})  to a later stage (Lemma 3).


We now show the self majorance property of  $ p(z) $.  It will be enough to verify  
that all coefficients in the expansion of  $ p(z) $  around zero are non negative.
Since  $  4\lambda( 1 - \lambda -\mu ) / ( 1-\mu ) $,  $ 2\lambda $, 
$ \mu\pi^2 / 4 $  and  $ \mu(1 - \mu) $  are all non negative, it is sufficient 
to prove that the Taylor coefficients of  $ L(x) $  are non negative.
For this aim we recall the formula  \cite[Section 3.14]{Spanier} 
%
%
\begin{equation}				\label{eq:tan-lambda}	
 \frac{ \tan \pi x / 2 } { \pi x /2 } 
	= \frac{8} {\pi^2} \sum_{k=0}^\infty \lambda(2k+2) x^{2k}, 
	   \qquad	\lambda(p) = \sum_{n=0}^\infty  \frac{1} { (2n+1)^p }  \  . 
\end{equation} 
According to  (\ref{eq:tan-lambda}), 
\begin{equation}				\label{eq:L(X)-expansion}  
    L(x) = \left( \frac{\pi }{ 2 } \right) \tan \left( \frac{\pi x }{ 2 } \right) 
		- \frac{2\lambda }{1-\mu} \, \frac{x}{1 - x^2 } 
=  \sum_{k=0}^\infty \Big( 2\lambda(2k+2) - \frac{2\lambda }{1-\mu} \Big) x^{2k+1} . 
\end{equation}
Since  $ \lambda(p) > 1 $, the coefficients of the power series of  $ L(x) $  are 
\begin{equation*}
  2 \lambda(2k+2) - \frac{2\lambda }{1-\mu} > 2 - \frac{2\lambda }{1-\mu} 
	= \frac{ 2(1 - \mu - \lambda) }{1-\mu}  \ge  0 .
\end{equation*} 
Consequently  $ p(z) $  is self majorant.

The univalence of  (\ref{eq:4.27})  follows as in Theorem 1. 
Consider now the issue of sharpness of  (\ref{eq:4.26}).   Note that  
$ \lambda + \mu $  is the multiplicity of the zeros of   
    $ u(x)  = (1 - x^2)^\lambda \cos^\mu(\pi x /2) $  at  $ x=1, -1 $.
From  (\ref{eq:4.34})  one easily deduce that $ \lambda \ge (1 - \mu)/2 $ 
which implies  $ \lambda + \mu \ge (1 + \mu)/2  $.  
Since we assumed that  $ \mu >0 $,  it follows that 
\begin{equation}				\label{eq:lambda+mu>1/2} 
	 1/2 \le \lambda + \mu \le 1 . 
\end{equation}
Consequently, 
$$
\int^1 \frac{ dx }{ (1 - x^2 )^{ 2\lambda } \cos^{ 2\mu } (\pi x / 2) } = \infty . $$
Hence the sharpness of  (\ref{eq:4.26})  follows.
Thus the proof of the theorem is complete, except of inequality  (\ref{eq:L'-ineq}).
\qed

\noindent
{\bf Remark.}  \ 
We saw that  $ \lambda < 0 $  is not possible.  
Due to symmetry between  $ \cos (\pi x / 2) $   and  $ 1 - x^2 $  it is natural 
to think that  $ \mu < 0 $  cannot occur as well.  But this is not true. 
It is easily seen that if  $ \lambda > 1/2 $  and  $ \mu < 0 $   is small enough 
in absolute value then  $ p(z) $  is self majorant  and  $ \varphi' \le 0 $.  
Thus our theorem is extended for the case   $ \lambda > 1/2 , \mu < 0 $  
and small in absolute value.  Moreover, it is possible to find  
concrete estimates for the case  $ \mu < 0 $,  but the somewhat involved details 
will be omitted.

Finally we mention that the particular case  $ \lambda > 0 $,  $ \mu >  0 $,
$ \lambda + \mu =1 $  was mentioned by Beesack  \cite{Beesack}  but he did not 
provide the proof.

The next lemma provides some inequalities about the function  $ G(x) $  
of  (\ref{eq:G=(1-x^2)tan})  which are needed to prove  (\ref{eq:L'-ineq}).

\begin{lemma}		
					\label{lemma:3} 
For the interval (0,1) we have:	 \\
(i)   \	The function  $  \  G(x) := ( 1 - x^ 2 ) \tan( \pi x / 2)  \  $  
of  (\ref{eq:G=(1-x^2)tan})    is convex.	\\
(ii) \ The function 
	$  \  \dfrac { ( 1 - x^ 2 ) \tan( \pi x / 2) } { \pi x / 2 }  \  $ 
decreases for  $ 0 \le x \le 1 $  and its maximal value at  $ x=0 $  is 1.	\\
(iii)
$$	 G(x)  \le  \frac{2}{\pi} (1 + x) . $$ 
(iv) 
$$	 \frac{\pi}{2} G'(x)   \le  1 + \frac{\pi^2}{4} (1 - x^2 ) . $$
\end{lemma}

\noindent
{\bf Proof. }    \ 
Convexity of  $ G(x) $  in the interval (0,1) follows from  (\ref{eq:tan-lambda}). 
Indeed, 
\begin{equation} 				\label{eq:G-series}	
 G(x) = \frac{ \pi } { 2 } x  
   - \frac{4} {\pi} \sum_{k=1}^\infty  \left( \lambda(2k) - \lambda(2k+2) \right)  x^{2k+1}
\end{equation} 
and the second derivative is negative due to monotonicity of  $ \lambda(p) $. 
This proves (i).

Since  $ \left(  \dfrac{ G(x) } { \pi x / 2 }  \right)'  \le  0 $  by (\ref{eq:G-series}), 
it follows that   $ \dfrac{ G(x) } { \pi x / 2 }  $  is decreasing and thus  its 
maximum is attained at  $ x=0 $.  This proves  (ii).

For the proof of  (iii)   we recall that   $ G(1) = 4/\pi $,  $ G'(1) = 2/\pi $  
as calculated in  (\ref{eq:G(1)}).  The equation of the tangent line through the 
point  $ ( 1  , 4/\pi ) $  is   
     $ y = \frac{4}{\pi} - \frac{2}{\pi} \, (1-x) = \frac{2}{\pi} \, (1+x) $. 
Now (iii) follows using the convexity proved in  (i).

In order to establish  (iv)  we use  (iii)  twice. Indeed,
$$
    G'(x) = -2x \tan( \pi x / 2)
        + \frac{\pi}{2} (1 - x^2) \left( 1 + \tan^2( \pi x / 2) \right) ,
$$
hence
\begin{align*}
G'(x) - \frac{\pi}{2} (1 - x^2) & = -2x \tan( \pi x / 2)
                + \tan( \pi x / 2) \frac{\pi}{2} G(x)    \\
&  \le      -2x \tan( \pi x / 2) + (1+x)  \tan( \pi x / 2)        \\
& = \frac {1 - x^2 }{1 + x} \tan( \pi x / 2) =  \frac{G(x)}{1 + x}
 \le \frac{2}{\pi} .
\end{align*}
This ends the proof of  (iv).	\qed

Now we are ready to complete the proof of  (\ref{eq:L'-ineq}).  First recall that
$ L(x) = (\pi/2) \tan ( \pi x / 2 ) 
	- \frac{2\lambda }{1-\mu} \, \frac{x}{1 - x^2 } \ge 0 $
for  $ 0 \le x \le 1 $.
This follows from the fact that the coefficients of the power series for  $ L(x) $
in  (\ref{eq:L(X)-expansion})  are positive.   Second,
$ \dfrac{-2\lambda} {1 - \mu} \le -1 $  by  (\ref{eq:4.34}).
In addition, by  (ii)  we have
$$
\frac{\pi}{2} (1 - x^2) \tan \left( \frac{\pi x} { 2 } \right)
= \left( \frac{\pi}{2} \right)^2  \dfrac { ( 1 - x^ 2 ) \tan( \pi x / 2) } { \pi x / 2 } \, x
\le      \left( \frac{\pi}{2} \right)^2 x .
$$ 
Moreover, by (iv), 
$$
\frac{\pi}{2} \left( (1 - x^2)  \tan \frac{\pi x} { 2 } \right)' 
	=  \left( \frac{\pi}{2} \right) G'(x)   \le   1 + \frac{\pi^2}{4} (1 - x^2 ) .
$$
Combining these  facts together,  we finally have
\begin{align*}
2  &  \left[ ( 1 - x^2 ) L(x) \right] \left[ ( 1 - x^2 ) L(x) \right]'	\\
&   \le  2 \left(  \frac{\pi}{2} (1 - x^2)\tan \left( \frac{\pi x} { 2 } \right) 
		- \frac{2\lambda }{ 1 - \mu } x \right)
\left(  \frac{\pi}{2} \left( (1 - x^2)\tan \left( \frac{\pi x} { 2 } \right) \right)' 
		- \frac{2\lambda }{ 1 - \mu } \right)		\\
&  \le  2 \left(  \left( \frac{\pi}{2} \right)^2 x - x \right)
    \left( 1 + \frac{\pi^2}{4} (1 - x^2 ) - 1 \right)		\\
& =   2 \left(  \left( \frac{\pi}{2} \right)^2  - 1 \right) \frac{ \pi^2}{4}  x (1 - x^2) 
	\le \pi^2 x (1 - x^2)  .
\end{align*}
Since  $ 0 < \mu \le 1 $,  the last inequality implies  (\ref{eq:L'-ineq}).

The third theorem corresponds to case (c) which is generated by 
the function  $ u(x) = (1 - x^2)^a / (1 + x^2)^{b} $.


\begin{theorem}
(A)  \  Let
\begin{equation}                	\label{eq:p-example3}    
 p(x) = \frac{  (2a - 4ab) + x^2 (2a+4ab-4a^2) } { (1 - x^2)^2 } 
	+ \frac{ (2b+4ab) + x^2 (-2b-4b^2 +4ab) }  { (1 + x^2)^2 } 
\end{equation}
and let   $ a, b $   satisfy
\begin{equation}                	\label{eq:example3-condition}  
   \frac{1}{2} \le a \le 1,	\qquad 
  \frac{1}{4}\left( -(5+4a) + (25+48a)^{1/2} \right) \le b \le a - \frac{1}{2} .
\end{equation}
  Then if  $ f(z) $  is an analytic function in   $ \Delta $  satisfying
\begin{equation}                	\label{eq:thm3-bound}        
    | Sf(z) | \le 2 p( |z| ) , \qquad  z \in \Delta ,
\end{equation}
it follows that  $ f(z) $  is univalent in  $ \Delta $.        \\
(B)  \  Let  $ p(x) $  be as in  (\ref{eq:p-example3})  and we have the 
more restrictive conditions 
\begin{equation}                	\label{eq:4.48}        
\begin{aligned}
   &   \frac{1}{2} \le a \le 1,			\\ 
  & \max \Big\{ a-1,  \frac{1}{4}\left( -(5+4a) + (25+48a)^{1/2} \right) \Big\} 
   \le   b 
   \le   \min \Big\{ a - \frac{1}{2},  
		- \frac{1}{2} + \frac{1}{2} \Big( 1 +4a(1-a) \Big)^{1/2} \Big\} ,
\end{aligned}
\end{equation}
then  $ p(z) $  is self majorant and
\begin{equation}                	\label{eq:4.49}        
   f(z) = \int_0^z \frac{dt} { u^2(t) }
        = \int_0^z  (1 - t^2)^{-2a}  (1 + t^2)^{2b}    \, dt
\end{equation}
is an odd univalent function in  $ \Delta $.  Moreover, the condition
 (\ref{eq:thm3-bound})  is sharp.
\end{theorem}

\noindent
{\bf Proof.} \ 
We  use the same reasoning as in Theorem 1 and 2.  Hence we start with  $ p(x) $ which is
given in (\ref{eq:p-example3}) and look for conditions ensuring that  $ p(x)(1 - x^2)^2 $  
is positive and decreases for  $ 0 \le x  \le 1 $.   It will be convenient to denote
	$  \ t = (1-x^2)/(1+x^2)  \ $ 
for  $ 0 \le x \le 1  $. 
With this notation and denoting  $ \psi(t) = p\big( x(t) \big) ( 1 - x^2(t) )^2 $, 
$ x^2(t) = (1-t)/(1+t) $  for  $ 0 \le t \le 1 $,  we have 
\begin{equation*}
\begin{aligned}	
  \psi(t) =  &   (2a - 4ab) + (1-t)(1+t)^{-1}(2a+4ab-4a^2) 	\\
	     &	+ t^2\left( (2b+4ab) + (1-t)(1+t)^{-1} (-2b-4b^2 +4ab) \right) .
\end{aligned} 
\end{equation*}
$ p(x)(1 - x^2)^2 $  decreases if  $ \psi'(t) \ge 0 $  since  $ dt/dx  \le 0 $ 
for  $ 0 \le x  \le 1 $.

We first consider the four necessary conditions
\begin{equation}				\label{eq:thm3-p(0)}	
   p(0) \ge 0,  \quad  p(1) \ge 0,  \quad  \psi'(0) \ge 0,  \quad  \psi'(1) \ge 0. 
\end{equation}
$ p(0) \ge 0 $  leads to 
\begin{equation}				\label{eq:thm3-a+b}	
		a+b \ge 0 .
\end{equation} 
From   $ p(1) \ge 0 $  and  $ \lim_{x \to 1} p(x)(1 - x^2)^2 =4a(1-a) $ 
it follows that 
\begin{equation}				\label{eq:thm3-a}	
		0 \le a \le 1 .
\end{equation} 
We have
\begin{equation*}
\begin{aligned}
  \psi'(t) =  &   -2(2a+4ab-4a^2) (1+t)^{-2}	\\
	& +2t \left( (2b+4ab) + (1-t)(1+t)^{-1} (-2b-4b^2 +4ab) \right)	 \\
	& -2t^2 (-2 b - 4 b^2 + 4ab) (1+t)^{-2} .
\end{aligned} 
\end{equation*}
$ \psi'(0) \ge 0 $  implies
\begin{equation}				\label{eq:thm3-a-ge-b}	
	a \ge b +  1/2
\end{equation} 
and from  $ \psi'(1) \ge 0 $  we get
	$ \ 5b -a + 4ab + 2(a^2 + b^2) \ge 0 $. 
The roots of the quadratic equation  $ \  5b -a + 4ab + 2(a^2 + b^2) = 0  \ $   
are
  $  \ b_{1,2} = \left( -(5+4a) \pm  (25+48a)^{1/2} \right) / 4  \  $  
and  
we observe that  $ b \le b_1 $  is impossible due to  (\ref{eq:thm3-a+b}). 
So we are left with the condition  $ b \ge b_2 $,  i.e., 
\begin{equation}				\label{eq:thm3-b-ge}	
	b  \ge  \frac{1}{4} \left( -(5+4a) + (25+48a)^{1/2} \right) .
\end{equation} 
Note that  (\ref{eq:thm3-b-ge})  includes  (\ref{eq:thm3-a+b}).

We now show that necessarily   $ a \ge 1/2 $.   Indeed, from  (\ref{eq:thm3-a-ge-b}) 
and  (\ref{eq:thm3-b-ge}),
\begin{equation*}
  a  \ge  b + \frac{1}{2}  
     \ge  \frac{1}{2} + \frac{1}{4} \left( -(5+4a) + (25+48a)^{1/2} \right)
\end{equation*}
or  $ 8a+3  \ge  (25+48a)^{1/2} $,   which implies  $ a \ge 1/2 $.

We conclude the following  necessary conditions
\begin{equation*}
  b \ge a - \frac{1}{2},  \quad   \frac{1}{2} \le a \le 1,	\quad 
  b \ge \frac{1}{4} \left( -(5+4a) + (25+48a)^{1/2} \right) . 
\end{equation*}
These are exactly conditions (\ref{eq:example3-condition})  of Theorem 3.  
See Figure 3.

It is our aim now to show that these condition are, indeed, sufficient to imply
monotonicity of  $ \psi(t) = p(x)(1 - x^2)^2 $  and thus proving part (A) of Theorem 3.
For this aim it will be convenient to verify the positivity of 
	$ H(t) := \frac{1}{4} (1 + t)^2 \psi'(t) $.  
 Using  the (trivial) relation  $ t^2 +(1-t^2) = 1 $,   $ H(t) $  is written as 
\begin{equation}
\begin{aligned}					\label{eq:thm3-H}
    H(t) 
= & \frac{1}{4} (1 + t)^2\psi'(t) 	\\
= & \left[ (1 - t^2) (2a^2 - 2ab - a) + t^2 (2a^2 - 2ab - a) \right] + t^2 (b + 2b^2 - 2ab)  
		\\
  &  \quad   + t(1 + t)^2 (b + 2ab) + t(1-t^2) (-b-2b^2 +2ab) 	\\
= &  (1 - t^2) (2a)(a-b-1/2) + t^2 (a-b)(a-b-1/2) 	\\
  &  \quad   +  t(1 + t)^2 (b + 2ab) + t(1-t^2) (2b) (a-b-1/2) .
\end{aligned} 
\end{equation}
The first two terms of  (\ref{eq:thm3-H})  are positive for  $ 0 \le t \le 1 $
since    $ (2a)(a-b-1/2) \ge 0 $   and   
$ (2a^2-2ab-a) + (b + 2b^2 - 2ab) =  (a-b)(a-b-1/2) \ge 0 $  by  (\ref{eq:thm3-a-ge-b}).

Now we consider the other two remaining terms.   If  $ b \ge 0 $  then  $ 2(b+2ab) $ 
 is obviously positive and  $ 2b (a-b-1/2) $  is positive by  (\ref{eq:thm3-a-ge-b}).
This ends the case  $ b \ge 0 $.

Assume now  $ b \le 0 $.  Here we use a different idea to show the positivity 
of  $  \psi'(t) $  (which we are after) or, equivalently, the positivity of 
  $ H(t) = \frac{1}{4} (1 + t)^2 \psi'(t) $.   Since   
$$ 
    H(0) = \frac{1}{4} \psi'(0) \ge 0,  \qquad   H(1) = \psi'(1) \ge 0 ,
$$ 
it is sufficient to confirm that  $ H(t) $  is convex,  i.e.,  $ H''(t) \le 0 $
on the interval  $(0,1) $.  Indeed, we have from  (\ref{eq:thm3-H}) 
\begin{equation*}
  	H''(t) = 12t(b+b^2) + 2b(2a+2b+3) = 2b \Big[ 6t (b+1) + 2(a+b) = 3 \Big].
\end{equation*}
But  $ b \le 0 $,  $ b+1 \ge b+a \ge 0 $  (since  $ a+b \ge 0 $  by  (\ref{eq:thm3-a+b})), 
and thus the proof that  $ H''(t) \le 0 $   is complete.  This ends the proof of part (A).

We now turn to the proof of part (B).  Hence we look for conditions ensuring 
positivity of coefficients of  $ p(x) $.  From  (\ref{eq:p-example3})  we get 
using  $ (1-t)^{-2} = \sum_{n=0}^\infty (n+1) t^n $,  
$ (1+t)^{-2} = \sum_{n=0}^\infty (n+1) (-t)^n $,  that
\begin{equation*} 			\label{eq:thm3-p(t)}	
\begin{aligned}
\frac{1}{2} p(x)  
	 = & \left( (a-2ab) + t(a+2ab-2a^2) \right) \sum_{n=0}^\infty (n+1)t^n   \\
	   & + \left( (b+2ab) + t(-b-2b^2+2ab) \right) \sum_{n=0}^\infty  (n+1)(-t)^n  . 
\end{aligned}
\end{equation*}
We already know that  $ p(x) \ge 0 $.  Thus the free coefficient is positive.  
For the coefficient of  $t$  we get the condition 
\begin{equation}			\label{eq:thm3-coeff-1}  
	3(a-b)) \ge 2(a+b)^2.
\end{equation}
In addition
\begin{equation}			\label{eq:thm3-coeff-2m}  
	(a+b)(1-a+b) \ge 0,
\end{equation}
\begin{equation}			\label{eq:thm3-coeff-2m+1}  
	a-b \ge a^2 + b^2 ,
\end{equation}
for the coefficients of  $ t^{2m} $  and  $ t^{2m+1} $,  respectively.

So our aim now is to find the additional constrains on  $ a,b $  so that 
(\ref{eq:thm3-coeff-1}), (\ref{eq:thm3-coeff-2m}),  (\ref{eq:thm3-coeff-2m+1})
will be satisfied.  For this we need to solve the two quadratic equations arising 
 naturally from  (\ref{eq:thm3-coeff-1}), (\ref{eq:thm3-coeff-2m+1}) 
and take into account also the condition  $ 1-a+b \ge 0 $  
(which is equivalent to  (\ref{eq:thm3-coeff-2m})  due to  $ a+b \ge 0 $).

We leave these straight forward computations to the reader. Summing up and 
recalling the previous restrictions from part (A) we finally get for the lower 
and upper bounds
\begin{equation}			\label{eq:thm3-MAX-1}  
	 b \ge  \max \big\{ a-1,  \  \frac{1}{4} (-5-4a+(25+48a)^{1/2} )  \big\},
\end{equation}
\begin{equation}			\label{eq:thm3-MAX-2}  
	 b \le  \min \big\{ a-\frac{1}{2},  \  \frac{1}{2} (-1+(1+4a-4a^2)^{1/2} ), 
		            \     \frac{1}{4} (-4a-3+(9+48a)^{1/2} )  \big\},
\end{equation}
Next it is needed to show that all three competing bounds appearing in 
(\ref{eq:thm3-MAX-2})   
are equal at the point  $  a = a_0 = (\sqrt{3}+1)/4 $ = 0.6830.  
Indeed - this somewhat surprising fact - follows easily by comparing the bound  
$ a-1/2 $  with each of the two other bounds (solving two quadratic equations, etc.)
Also it follows from the same computation that  $ a-1/2 $  is the minimum among 
the three bounds for the  interval  $ [1/2, a_0] $.
From what was said above for the interval  $ [1/2, a_0] $  it follows that  
$ a-1/2 $  is smaller than the two other bounds.
This ends the discussion for the interval  $ [1/2, a_0] $.

Now consider the other interval $ [a_0,1 ] $.
In this interval  $ a-1/2 $   is obviously not competing any more since  $ a-1/2 $
 is greater than the the two other bounds.  Hence it is left to decide which bound 
between the two other bounds is the ``winner" for the remaining interval  $ [a_0,1 ] $. 
Denote the two bounds by
\begin{equation*}			\label{eq:thm3-MAX-Ta-Ga}  
  T(a) = (1/4)(-4a-3+(9+48a)^{1/2} ), \qquad  S(a) = (1/2)(-1+(1+4a-4a^2)^{1/2} )  .
\end{equation*}
From our earlier considerations  $ T(a_0) = S(a_0) $.  

We now show that  $ S(a) \le  T(a) $  in the interval  $ [a_0,1 ] $,  i.e.,  $ S(a) $
 ``wins" for the interval $ [a_0,1 ] $.  
$ S(a) \le  T(a) $  is equivalent to  
	$ 1 + 4a +2(1 + 4a - 4a^2)^{1/2} \le  (9+48a)^{1/2} $, 
squaring this implies 
	$ (1 + 4a) (1 + 4a - 4a^2)^{1/2} \le  1 + 6a $, 
i.e., 
\begin{equation}			\label{eq:thm3-4.72}  
   	1 + 4a - 4a^2  \le  \left( \frac{ 1+6a }{1+4a} \right)^2 . 
\end{equation}
  It is immediate to check that the derivative 
of the left hand side and the right hand side of  (\ref{eq:thm3-4.72})  
are negative and positive respectively.
 Since  $ T(a_0) = S(a_0) $,  it now follows that   $ S(a) \le T(a) $  for 
the interval   $ [a_0,1 ] $.  This ends the discussion for the upper bound
(\ref{eq:thm3-MAX-2}).

The treatment of the lower bound  (\ref{eq:thm3-MAX-1})  case is much shorter. 
Indeed, from  (\ref{eq:thm3-MAX-1})   
we need to equate only two bounds and not three as in the case of the upper bound.
It is easily seen that for  $ a_1  = (\sqrt{7} + 1)/4 = 0.9114 $  we get equality 
between the two bounds. Thus we easily conclude that
\begin{equation*}
\begin{aligned}
	(1/4)(-5-4a+(25+48a)^{1/2} ) \le  b  \qquad  &  \hbox{\rm if} 
			\qquad 1/2 \le a \le ( \sqrt{7} + 1 )/4	,	\\
	a-1 \le b \qquad 	& \hbox{\rm if} 
			\qquad ( \sqrt{7} + 1 )/4  \le a \le 1 .
\end{aligned}
\end{equation*}
Hence the major part of the proof of part (B) is complete.  See Figure 4.

The univalence of  $  \displaystyle    f(z) = \int_0^z \frac{dt}{ u^2(t) } 
   = \int_0^z  \frac{ dt } { (1 - t^2)^{2a} (1 + t^2)^{-2b}  }  \  $
follows by similar reasoning as in Theorems 1 and 2.  The sharpness of 
(\ref{eq:thm3-bound})  follows from  $ a \ge 1/2 $,  which implies divergence 
of the integral at  $ z=1 $.
\qed

Finally we remark that Nehari considered partial case of part (A) in
\cite{Nehari-1979}.   He used the letters  $ \mu, \lambda $.  The connection 
with our notation is  $ a = (1 + \mu) / 2, \ b= \eta /2 $. 
  Nehari considered the special two cases  $ \mu = \eta $ and  $ \eta=0 $, 
or in terms of our notation  $ b=0 $  and   $ a = b+ 1/2 $.

\bigskip
\begin{center}
    \includegraphics[scale=0.8]{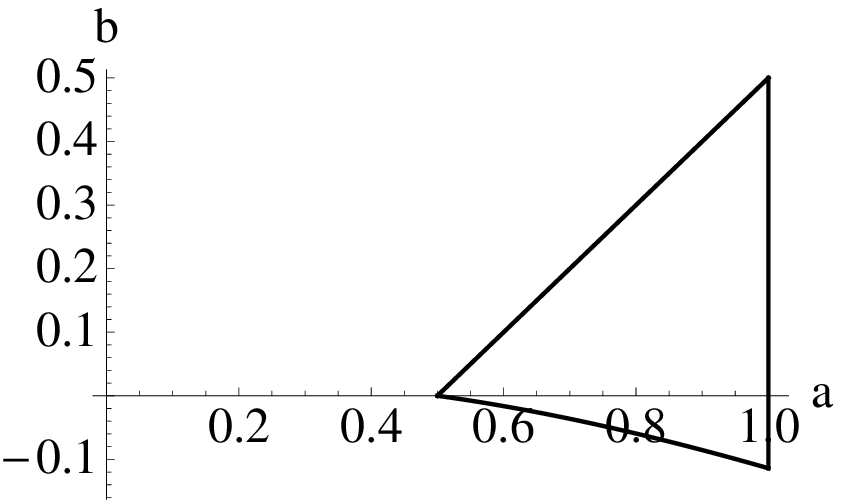}
\end{center}
\centerline{ Figure 3.  The domain restricted by  (\ref{eq:example3-condition})}

\bigskip
\begin{center}
    \includegraphics[scale=0.8]{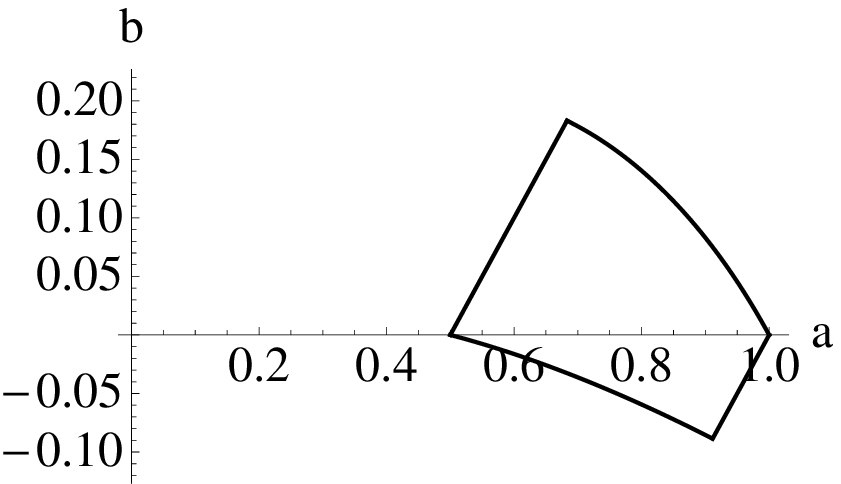}
\end{center}
\centerline{ Figure 4.  The domain restricted by  (\ref{eq:4.48}) }


\section{ Two sufficient conditions for univalence arising 	\\
		from Lemma 2. }

We start with the condition arising from case  (d),  which is generated by 
the function  $ u(x) = \cos(\pi x/2) e^{ \lambda x^2 } $.


\begin{theorem}
  Let
\begin{equation}                	\label{eq:p-example4}    
   p(x) =  2\lambda \pi x \tan(\pi x /2) + \pi^2 / 4 - (2\lambda + 4\lambda^2 x^2)
\end{equation}
and let   $ \lambda $   satisfy
\begin{equation}                	\label{eq:example4-condition}  
  0  \le \lambda  \le \lambda_0 = \frac{1}{8}\left( (4+\pi^2) - (16+\pi^4)^{1/2} \right). 
\end{equation}
  Then if  $ f(z) $  is an analytic function in   $ \Delta $  satisfying
\begin{equation}                	\label{eq:5.3}        
    | Sf(z) | \le 2 p( |z| ) , \qquad  z \in \Delta ,
\end{equation}
then  $ f(z) $  is univalent in  $ \Delta $  and condition  (\ref{eq:5.3}) is sharp.        
\end{theorem}

\noindent
{\bf  Proof. } 
We first point out that  $  \lambda \ge 0 $  is necessary. 
Indeed, this follows at once from the  condition  $ p(x) \ge 0 $  
by letting   $ x \to 1 $  and  noting that  $ \pi x \tan(\pi x /2) \to \infty $ 
while the other terms remain finite.

Next, for  $ \varphi = (1-x^2)^2 p(x) $   we have by direct calculation 
that 
\begin{equation*}
\begin{aligned}
  \varphi'(x) 
&  = \tan^2(\pi x /2) \Big( \lambda \pi^2 x (1-x^2)^2 \Big) 
   + \tan(\pi x /2) \Big( 2\lambda \pi (1-x^2)^2 -4x (1-x^2) (2\lambda \pi x) \Big) \\
&  + x\left( (1-x^2)^2 (\lambda \pi^2 - 8 \lambda^2  \right) 
   - \pi^2 x (1-x^2)  +  8\lambda x (1-x^2) + 16\lambda^2 x^3 (1-x^2) .
\end{aligned}  
\end{equation*}
For later use we summarize that
\begin{equation}                	\label{eq:example4-phi'/x}   
\begin{aligned}
  \frac{ \varphi'(x) } {x} 
&  = \left(  \frac{ (1-x^2) \tan(\pi x /2) } { \pi x /2 }  \right)^2
     	 \lambda \pi^2  (\pi/2)^2  x^2 		\\
&   + \left(  \frac{ (1-x^2) \tan(\pi x /2) } { \pi x /2 }  \right) 
     \left( 2\lambda \pi (1-x^2) \left( \frac{\pi}{2} \right) - 4\lambda \pi^2 x^2 \right)  
	\\
&  + \left[  (1-x^2)^2 ( \lambda \pi^2 - 8 \lambda^2 ) 
   - \pi^2  (1-x^2)  +  8\lambda  (1-x^2) + 16\lambda^2 x^2 (1-x^2) \right].
\end{aligned}  
\end{equation}
Since  $ \varphi'(x) \le 0 $   is assumed in Theorem A,  we have
from  $ \varphi'(x) / x \le 0 $  as  $ x \to 0^+ $, 
\begin{equation}			\label{eq:example4-quadratic}   
	2 \lambda \pi^2 - 8 \lambda^2 - \pi^2 + 8 \lambda   \le   0 .
\end{equation}
The roots of 
         $  \ 2 \lambda \pi^2 - 8 \lambda^2 - \pi^2 + 8 \lambda = 0  \ $ 
are 
         $  \ \lambda_{0,1} = \left( (4+\pi^2) \pm (16+\pi^4)^{1/2} \right) / 8  \ $. 
$ \lambda \ge \lambda_1 $   is impossible in view of  $ p(0) \ge 0 $, 
i.e.,  $ \lambda \le \pi^2 / 8 $ .   Thus we are left with the condition 
(\ref{eq:example4-condition}), i.e., 
\begin{equation}			\label{eq:example4-lambda}   
    \lambda \le \lambda_0 = \left( (4+\pi^2) - (16+\pi^4)^{1/2} \right) / 8 = 0.40235...
\end{equation}

Our next aim is to show that (\ref{eq:example4-lambda}) is not only necessary 
but also sufficient for monotonicity of  $ \varphi(x) $.

  With the notation  $ R(x) := \dfrac{ (1-x^2) \tan(\pi x /2) } { \pi x /2 } $,
(\ref{eq:example4-phi'/x})  reads 
\begin{equation}                	\label{eq:example4-phi'-R}   
\begin{aligned}
  \frac{ \varphi'(x) } {x} 
&  = \frac{ \lambda \pi^4  }{4} x^2 R^2(x) 
  +  \left[ \lambda \pi^2 (1-x^2) - 4\lambda \pi^2 x^2 \right] R(x)	\\
&  + (1 - x^2) \left[ \lambda \pi^2 - 8 \lambda^2 - \pi^2 +  8\lambda 
			 + x^2 ( -\lambda \pi^2 + 24\lambda^2 ) \right].
\end{aligned}  
\end{equation}
From Lemma 3(ii) we know that  $ R(x) $  is decreasing for  $ 0 \le x \le 1 $  
and as   $ R(1) = \dfrac{8}{ \pi^2 }, \ R(0) = 1 $,  we have  
	$ \dfrac{8}{ \pi^2 } \le  R(x)  \le  1 $.
%
Hence we apply  $ R^2(x) \le R(x) $  and  $ R(x) \ge {8} /{ \pi^2 } $  
for two of the terms appearing in  (\ref{eq:example4-phi'-R})  and obtain 
\begin{equation}
\begin{aligned}
  &   \lambda  (\pi^4 /4) x^2 R^2(x)   - 4\lambda \pi^2 x^2 R(x)  	\\
  &  \le  \big[ \lambda  (\pi^4 /4 ) x^2  - 4 \lambda \pi^2 x^2 \big]  R(x)
	  = -\lambda \pi^2 (4 - \pi^2 / 4) x^2 R(x) 		\\
  &  \le  - \lambda \pi^2 (4 - \pi^2 / 4) x^2  \,  \frac{8}{ \pi^2 } 
    	  = -\lambda  ( 32 - 2 \pi^2 ) x^2  
\end{aligned} 
\end{equation}
which is negative;  For another term of  (\ref{eq:example4-phi'-R})  we apply  
$ \lambda \pi^2 (1-x^2) R(x) \le \lambda \pi^2 (1-x^2) $. 
Thus, from  (\ref{eq:example4-phi'-R})  and the above,  also using 
$ -(1-x^2)^{-1} \le -1 $, we have 
\begin{equation}			\label{eq:example4-phi'-end}   
\begin{aligned}
 \frac{ \varphi'(x) } { x(1-x^2) }   & 
    \le  
   -\lambda ( 32 - 2 \pi^2 ) x^2  + \lambda \pi^2
   +\left[ \lambda \pi^2 - 8 \lambda^2 - \pi^2 +  8\lambda 
			 + x^2 ( -\lambda \pi^2 + 24\lambda^2 ) \right]		\\
&  = ( 2\lambda \pi^2 - 8 \lambda^2 - \pi^2  + 8\lambda ) + 
		\lambda x^2 ( -32 +  \pi^2 + 24 \lambda ) .
\end{aligned} 
\end{equation}
The constant term on the right hand side of  (\ref{eq:example4-phi'-end})  is negative 
by (\ref{eq:example4-quadratic})  and the coefficient of  $ x^2 $   in  
(\ref{eq:example4-phi'-end})   is negative by  (\ref{eq:example4-condition}).  
Together, we conclude that   $ \varphi'(x) \le 0 $  and the proof of monotonicity 
is complete.

We now turn to prove that  $ p(z) $,  as presented in  (\ref{eq:p-example4}), 
is self majorant.  This is very easy.  Indeed, 
$ \lambda \le \lambda_0 \le 1/2 $   implies that the free coefficient in
the Taylor series of  $ p(x) $,  i.e.,  $ (\pi/2)^2 - 2\lambda $  is positive. 
It also implies the positivity of the coefficient of  $ x^2 $,  i.e.,  
$ \lambda( \pi^2 -4 \lambda) $.  
Positivity of all other coefficient follows at once from the  positivity of 
the Taylor coefficients of  $ \tan( \pi x/2) $. 
The proof of the remaining claims is identical with the proof of earlier 
theorems.   This ends the proof of Theorem 4. 
\qed

We start with the condition arising from (e), corresponding to the solution 
$ u(x) = \cos(\pi x /2) \exp( -\lambda \cos(\pi x /2) ) $.

\begin{theorem}
  Let
\begin{equation}                	\label{eq:p-example5}    
	p(x) =  \frac{\pi^2}{4}  
	\Big[ 1 - \lambda^2 \sin^2(\pi x /2) 
    		+ 2\lambda \frac{ \sin^2(\pi x /2) } { \cos(\pi x /2) } 
		- \lambda \cos(\pi x /2)  	
	\Big] .  
\end{equation}
and let   $ \lambda $   satisfy
\begin{equation}                	\label{eq:example5-condition}  
  0  \le \lambda  \le \lambda_0   = \frac{1}{\pi^2} 
  	\left[  ( 4 + \frac{5}{4} \pi^2 ) 
	   - \left( \Big( 4 + \frac{5}{4} \pi^2 \Big)^2  - 8 \pi^2 \right)^{1/2} 
 	\right] = 0.2664...
\end{equation}
  Then  $ p(z) $  is self majorant and if  $ f(z) $  is an analytic function 
in   $ \Delta $  satisfying
\begin{equation}                	\label{eq:example5-bound}        
    	| Sf(z) | \le 2 p( |z| ) , \qquad  z \in \Delta ,
\end{equation}
then  $ f(z) $  is univalent in  $ \Delta $  and the condition 
(\ref{eq:example5-bound}) is sharp.        
\end{theorem}


\noindent
{\bf Proof.}  
From  (\ref{eq:p-example5})  and the requirements  $ p(0) \ge 0 $, $ p(1) \ge 0 $,  
we easily get  $  1 -\lambda \ge 0 $, $ \lambda \ge 0 $.   But  $ \lambda = 1 $   is 
impossible, since it implies  $ p(x) \equiv 0 $.  Hence we get
	$ 0  \le   \lambda  < 1 $.
Note that  $ \lambda = 0 $   leads to the trivial case  $ u = \cos( \pi x/2) $ 
 and  $ p(x) = \pi^2/4 $.

It will be convenient to denote  $ \psi(x) = \dfrac{4}{ \pi^2 } \ p(x) (1-x^2)^2 $. 
Our aim is to find conditions ensuring monotonicity of  $ \psi(x) $,  i.e., 
$ \psi'(x) \le 0 $.   After some direct calculations we deduce from   
(\ref{eq:p-example5})  that
\begin{equation}                	\label{eq:example5-psi-prime}    
\begin{aligned}
  \frac{ \psi'(x) } { x(1-x^2) } = 
&    - 4 \left[ 1- \lambda^2 + \lambda^2 v^2 +2 \lambda/ v - 3\lambda v  \right] 	\\
&    + (3\lambda -2 \lambda^2 v) \ \frac{ \pi^2 }{4} \, 
		\frac{ \sin( \pi x/2 ) (1-x^2) } { \pi x/2 } 
	+ \lambda \ \frac{ \pi^2 }{2} \frac{1}{v}  \, 
		\frac{ \tan( \pi x/2 ) (1-x^2) } { \pi x/2 }  
\end{aligned}		     
\end{equation}
with  $ v(x) = \cos( \pi x /2 ) $.

First we check the condition  $ \psi'(0) \le 0 $.  From (\ref{eq:example5-psi-prime}) 
it follows that 
\begin{equation}                	\label{eq:example5-psi-prime-x=0}    
  \frac{ \psi'(x) } { x(1-x^2) }\Big|_{x=0, v=1}  = 
   -\frac{ \lambda^2 \pi^2 }{2} + \lambda\left( 4 + \frac{5}{4} \pi^2 \right) -4 \le 0 .
\end{equation}
The roots of the corresponding quadratic equation are 
\begin{equation*}                	\label{eq:example5-lambda-01}    
   \lambda_{0,1}  = \frac{1}{ \pi^2 } 
	\left[ 4 + \frac{5}{4}\, \pi^2  \pm
 	     \left( \Big(4 + \frac{5}{4}\, \pi^2 \Big)^2 - 8 \pi^2  \right)^{ 1/2 } 
        \right] .
\end{equation*}
$ \lambda \ge \lambda_1 $  is impossible because of  $ 0 \le \lambda < 1 $.  Hence
we are left with 
\begin{equation}                	\label{eq:example5-lambda-0}    
   \lambda  \le  \lambda_0  =  \frac{1}{ \pi^2 } 
	\left[ 4 + \frac{5}{4}\, \pi^2  -
 	     \left( \Big(4 + \frac{5}{4}\, \pi^2 \Big)^2 - 8 \pi^2  \right)^{ 1/2 } 
        \right] = 0.2664...
\end{equation}
appearing in  (\ref{eq:example5-condition}).

It is our aim now to show that (\ref{eq:example5-psi-prime-x=0})  
is not only necessary but also a sufficient condition for monotonicity.
For this aim note that 
	$ \sin( \pi x/2 ) / ( \pi x/2 ) \le 1 $ 
and from  Lemma 3, (ii),  it follows that
$ 	\dfrac{ \tan( \pi x/2 ) (1-x^2) } { \pi x/2 } 	\le  1 $.
Hence  (\ref{eq:example5-psi-prime})  and Lemma 3, (ii),  imply  
\begin{equation}                	\label{eq:example5-psi-prime-ineq}    
  \frac{ \psi'(x) } { x(1-x^2) }    \le 
    - 4 \left[ 1- \lambda^2 + \lambda^2 v^2 +2 \lambda/ v - 3\lambda v  \right] 	
    +  (3\lambda -2 \lambda^2 v) \ \frac{ \pi^2 }{4}  +  \frac{ \lambda \pi^2 }{2v} , 
\end{equation}
with equality for  $ x=0 $ (and  $ v=1 $).

In order to prove that  $ \psi' \le 0 $, it will be convenient to change the sign of 
the expression appearing on the right hand side of  (\ref{eq:example5-psi-prime-ineq}), 
consider it as a function of  $v$,  say  $ Q(v) $, and show that  $ Q(v) \ge 0 $  
for  $ 0 \le v \le 1 $.   Hence, in terms of  $ Q(v) $,  
(\ref{eq:example5-psi-prime-ineq})  is equivalent to 
\begin{equation}                	\label{eq:example5-Q-ineq}    
  Q(v) :=  4 - 4\lambda^2 + 4\lambda^2 v^2 
           + \frac{2 \lambda }{v} \left(  4 - \frac{ \pi^2 }{4}  \right)
	   - 12 \lambda v 	
    +  ( 2 \lambda^2 v - 3\lambda) \ \frac{ \pi^2 }{4}  \ge  0 .
\end{equation}
Note that inequality (\ref{eq:example5-psi-prime-x=0})  is identical with 
    $ Q(1) \ge 0 $.  
With the above investigation it follows that in order to verify that  $ Q(v) \ge 0 $ 
it is enough to show the easy inequalities  $ Q''(v) \ge 0,  \  Q'(1) < 0 $.
Indeed, this implies that  $ Q'(v) $  is negative and   $ Q(v) $  is decreasing.
Hence the positivity of  $ Q(1) $  implies  $ Q(v) \ge 0 $  and the proof
of the monotonicity is complete.

It is left to show that  $ p(z) $  in  (\ref{eq:p-example5})  is self majorant.
For that it is enough to check the positivity of the coefficients of  $ p(z) $ 
in  (\ref{eq:p-example5}).  
This is equivalent to the positivity of the Taylor coefficients of 
\begin{equation}                	\label{eq:example5-coeff}    
\frac{4}{\pi^2}  \,  p \Big( \frac{2}{\pi} x  \Big)    
	=  1 - \lambda^2  + 2\lambda \sec x 
	+ \frac{ \lambda^2 }{2} ( 1 + \cos 2x ) -3 \lambda \cos x .  
\end{equation}
We recall that 
\begin{equation}                	\label{eq:example5-sec-series}    
   \sec x = \sum_{ m=0 }^\infty  \frac{ | E_{2m} | }{ (2m)! }  x^{2m}, 
\end{equation}
where  $ E_{2m} $  are the Euler numbers $ E_0=1, \ E_2 = -1, \  E_4 = 5 $ 
and in general, 
\begin{equation}                	\label{eq:example5-euler-numbers}    
	E_{2m} = (-1)^m \left( \frac{2}{\pi} \right)^{2m+1}  (2m)!  \beta(2m+1),
\end{equation}
where 
$  \ \beta(2m+1) = \sum_{ n=0 }^\infty \dfrac{ (-1)^n } { (2n+1)^{ 2m+1 } }  \  $ 
is increasing to 1.

The free coefficient in the Taylor series of (\ref{eq:example5-coeff}),  $ 1 - \lambda $, 
is clearly positive.  The coefficient of  $ x^2 $,  i.e.,  
$  \lambda |E_2| -\lambda^2 + 3\lambda $  
is positive by $ |E_2| = 1 $  and the bound  (\ref{eq:example5-lambda-0})  on  $ \lambda $. 
For the coefficients of all other powers  $ x^{2m},  \  m \ge 2 $,  in (\ref{eq:example5-coeff}) we use 
\begin{equation}                	\label{eq:example5-coeff-ineq}    
\begin{aligned}
  \left\{  \frac{4}{\pi^2}  \,  p \Big( \frac{2}{\pi} x  \Big)  \right\}_{2m}   
	& =  \frac{ 1 }{ (2m)! } \Big( 2\lambda |E_{2m}| + \lambda^2 (-1)^m 2^{2m-1} - 3\lambda \Big) \\
	& \ge  \frac{ \lambda }{ (2m)! } \Big( 2 |E_{2m}| - \lambda 2^{2m-1} - 3 \Big) .
\end{aligned}
\end{equation}
We now apply  (\ref{eq:example5-euler-numbers})  to prove by induction that 
$ |E_{2m}| \ge 5^{m-1} $  with equality for  $ m=2 $.  The positivity of the required 
Taylor coefficients now follows from  (\ref{eq:example5-lambda-0})  and   
(\ref{eq:example5-coeff-ineq}).
This ends the proof of Theorem 5.
\qed

We end this chapter with an application.

The question of ``radius of univalence" is of some interest.  For example, some 
mathematicians were interested in finding the radius of univalence of the 
error function   
	$ \hbox{errf}(z) = \int_0^z e^{ -t^2} \, dt $. 
See  \cite{Kreyszig-Todd}  and more references there. 
By applying our criteria of univalence, we easily give a good (but not sharp) 
estimate for this radius.

In order to do that, our aim is to find an estimate for  $r$  such that 
	$ \hbox{errf}(rz) = \int_0^{rz} e^{ -t^2} \, dt $ 
is univalent for  $ |z| < 1 $.  We have 
		$  ( \hbox{errf}(rz) )' = r \exp( -r^2  z^2 ) $  
and thus for  $ Sf(z) = (f'' / f')'- \frac{1}{2}( f''/f')^2 $  \ 
we get
	$  S \big( \hbox{errf}(rz) \big) = -2 r^2( 1 + r^2 z^2 )  $.
The corresponding differential equation  $    u'' - r^2(1 + r^2 x^2 ) u = 0    $  
has solutions  
$ \,  y_1 = e^{ r^2 x^2 /2 }, \    y_2 = e^{ r^2 x^2 /2 } \int_0^{rx} e^{-t^2 } dt  \ $ 
with  $ y_2 / y_1 = \hbox{errf}(rx) $.  However  $ p_{errf}(x) = - r^2( 1 + r^2 x^2 ) < 0 $ 
is not a Nehari function, hence we have to compare  $ S \big( \hbox{errf}(rz) \big) $
 with another differential equation.

According to Theorem 5,  we get that a sufficient condition for univalence of
$ S \big( \hbox{errf}(rz) \big) $  is 
\begin{equation*}				\label{eq:errf-univ}  
	| S \big( \hbox{errf}(rz)  \big) | \le  2( r^2 + r^4 |z|^2 ) \le 2 p( |z| ) ,
\end{equation*}
where  $ p(x) $  is given in  (\ref{eq:p-example5})  with a suitable value of 
$ \lambda $. 
Hence, comparing the graph of the function   $ r^2 + r^4 x^2 $  (with  $ r=1.365 $) 
with the graph  of  $ p(x) $  (with  $ \lambda = 0.2 < \lambda_0 $),  we get for the 
radius  $r$  of univalence of the  errf$(z)$  the lower bound   $ r \ge 1.365 $. 
See Figure 5.   This is better than the various bounds for the radius of univalence 
which are quoted in  \cite{Kreyszig-Todd}  but, of course, it is not sharp. 
The precise radius of univalence is calculated by numerical methods in  \cite{Kreyszig-Todd} 
to be  $ r = 1.5748...$ \  .

\bigskip
\begin{center}
    \includegraphics[scale=0.8]{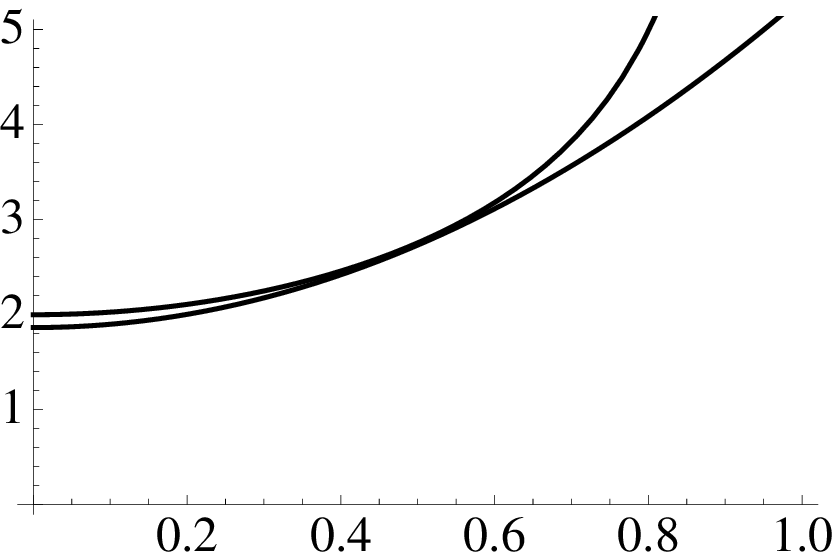}
\end{center}
\centerline{ Figure 5.    $ S \big( $errf$ (rx) \big) $   with  $ r=1.365 $ 
and  $ p(x) $  of  (\ref{eq:p-example5})   with  $ \lambda= 0.2 $ . }


\section{ Concluding remarks }

We end this paper with some observations and remarks.

First we recall the subject of quasiconformal extension.  As already mentioned in
the introduction, Ahlfors and Weill  \cite{Ahlfors-Weill}  were the pioneers of this 
subject. In their pioneering work the authors not only opened the new line of 
research concerning quasiconformal extension, but in fact gave an entirely new 
proof of Nehari's condition for univalence,  $ |Sf(z)| \le (1 - |z|^2)^{-2} $,
 based on ideas from the theory of  quasiconformal mappings. 
(See also the book of G. Shober   \cite[p. 169]{Schober}. 
Later on Gehring and Pommerenke wrote another deep paper on the subject.  
See  \cite[Thm. 4]{Gehring-Pommerenke}.
Among other things they extended the results of Ahlfors and Weil to any sufficient 
condition for univalence of the type  (\ref{eq:Sf<p}),  while Ahlfors and Weill considered 
the special case  $ p(z) = (1 - z^2)^{-2} $:  \\
{\sl 
If one has a sufficient condition for univalence such as
	$ | Sf(z) | \le 2 p( |z| ) $, 
then for any  $ 0 < a < 1 $,  the condition 
	$ | Sf(z) | \le 2a p( |z| ) $ 
is sufficient not only for univalence of  $ f(z) $  in the unit disc, but also for 
extension to a quasiconformal map in the full plane $ \cal C $.
}

This result may be obviously applied to all five theorems discussed in 
Sections 4 and 5.  
Further extensions of results of this nature were made by Chuaqui and Osgood  
\cite{Chuaqui-1995}, \cite{Chuaqui-1998}.

We now come back to the question of finite and infinite valence of an analytic 
function discussed in Section 2.
The remarkable example of Hille  mentioned in  (\ref{eq:Hille})  is the starting point 
of this subject.  As we already know, a function,  say  
$ f(z) =  \dfrac{ (1+z)^{ i\gamma } - (1-z)^{ i\gamma } } {(1-z)^{ i\gamma }} $,  
which satisfies 
\begin{equation}					\label{eq:6.3}
	 Sf(z) = 2 p(z),  \qquad  p(z) = ( 1 - z^2 )^{-2} (1 + \gamma^2 )  .
\end{equation}
has infinitely many zeros in the unit disc. 
A natural question arises:  If instead of condition (\ref{eq:6.3}),  a function  $ f(z) $  satisfies 
\begin{equation}					\label{eq:6.4}
	| Sf(z) |  \le 2 p( |z| ),  \qquad  p(z) = ( 1 - z^2 )^{-2} (1 + \sigma(z) )  , 
\end{equation}
for  $ |z| < 1 $,   where  $ | \sigma(z) |  \le \sigma( |z| ) $  and  $ \sigma(x) \searrow 0 $ 
as  $ x \to 1^- $,   is it still true that  $ f(z) $  may have infinitely many zeros in 
the unit disc?

In order to deal with this question we define for a Nehari's function  $ p(z) $
the quantity 
\begin{equation}					\label{eq:6.5}
	 \tau = \lim_{ x \to 1^- } p(x) (1 - x^2)^2 .
\end{equation}
In case  $ \tau < 1 $,  Gehring and Pommerenke proved  \cite[Thm. 4]{Gehring-Pommerenke} 
that a function  $ f(z) $   satisfying  $ | Sf(z) |  \le 2 p( |z| ) $ 
has finitely many zeros in the unit disc  $ \Delta $.  
On the other hand, if  $ \tau > 1 $  then a function  $ f(z) $  satisfying  
$  Sf(x) = 2 p(x) $   has infinitely  many real zeros. 
Indeed, if  $ \tau > 1 $, then  $ p(x) (1 - x^2)^2 \ge 1 + \gamma^2 $ 
on some interval  $ (r,1) $  with  $ r $  sufficiently close to 1 and for some  
$ \gamma > 0 $.
But Hille's equation  $ u'' + ( 1 - x^2 )^{-2} (1 + \gamma^2 ) u = 0 $  has the 
real solution  
  $ u(x) = (1 - x^2)^{1/2} \cos \left( \dfrac{\gamma}{2} \log \dfrac{1+x}{1-x} \right) $
with infinitely many real zeros in $ (r,1) $.
Now the claim follows easily from the Sturm comparison theorem.

It remains to discuss the case  $ \tau = 1 $. 
First we point out that if  $ \tau = 1 $,  the function  $f(z)$  can not be univalent unless
$ p(x) (1 - x^2)^2 \equiv 1 $.   See  \cite{Chuaqui-2007A}  for the proof which uses the 
so-called  ``relative convexity lemma''.  (Later on we reffer to the above as ``uniqueness 
argument").  Alternatively, one can replace the relative convexity lemma by the proof the 
Singular Sturm comparison theorem  \cite{Aharonov-Elias},  which is more direct.

We now turn to the question of finite or infinite valence for the case  $ \tau = 1 $
which needs more care.  A very nice contribution was made by Chuaqui, Duren, Osgood and Stowe 
 in  \cite{Chuaqui-2009}.  The authors showed that if the convergence to 1 is not ``too fast''
then it follows that the function is infinite valent, while if the solution 
is  ``fast enough''  then the situation is more delicate.  The authors restricted 
their analysis for the real case.  Theorem 1 in their paper yields at once the following result:  \\
{\sl 
Given Nehari's function  $ p(x) $  of the form 
\begin{equation*}
	p(x) = ( 1 - x^2 )^{-2} (1 + \sigma(x) )  , 	\qquad
	\sigma(x) = \lambda \left( \log \frac{1}{ 1 - x^2 }   \right)^{-2}, 
	\qquad 		0 < x < 1 .
\end{equation*}
If  $ \lambda > 1 $   then a function  $ f(x) $   satisfying 
	 $ Sf(x)  \ge 2 p(x) $, 
has infinitely many real zeros near  $ x=1 $.
On the other hand, if   $ \lambda < 1 $  and a function  $ f(x) $  satisfies 
	 $ Sf(x)  \le 2 p(x) $, 
then  $ f(x) $  has finite valence on  $ (0,1) $.
}

We now will examine what is the value of  $ \tau $  for each of our five theorems. 
First we consider Theorems 4 and 5.  In both cases  $ \tau $  is easily seen
to be zero due to the slow growth of  $ p(x) $  appearing in Theorems 4 and 5.
Obviously, multiplying  $ p(z) $  by any positive constant  $C$,  $ \tau $  is again zero. 
Thus by the result of Gehring and Pommerenke quoted above, the functions  $ f $  for which
\begin{equation}					\label{eq:6.11}
	| Sf(z) | \le 2C p( |z| ), \qquad 	C \quad  \hbox{any positive constant}, 
\end{equation}
are of finite valence.

But in fact, we can be more specific and give some quantitative result.
For this aim we use the following  result of Chuaqui, Duren and Osgood 
\cite[Theorem 2]{Chuaqui-2008}:

\noindent
{\bf  Theorem C.} 
{\sl 
Let  $f$  be analytic and locally univalent in  $ \Delta $,  and suppose its Schwarzian 
derivative satisfies
\begin{equation}					\label{eq:6.12}
	| Sf(z) | \le \frac{ 2C } { 1 - |z|^2 } , 
\end{equation}
for a constant  $ C > 2 $.  Then  $ f(z) $  has finite valence  $ N(C) < AC \log C $ 
where  $A$   is an absolute constant.	}

Now, for  $p(x)$  appearing in either Theorem 4 or Theorem 5,  condition (\ref{eq:6.12}) 
is fulfilled, as  $ \tan( \pi x/2  ), 1/ \cos( \pi x/2 ) $  are both bounded 
by  $ A/(1-x^2) $  for some absolute constant  $A$.  Hence the following statement follows:

{\sl 
If  $ p(x) $  is the function appearing in either Theorem 4 or Theorem 5,
we have for  $ P(x) = C p(x), C>2 $  and a function  $ f(z) $   satisfying
 $ | Sf(z) | \le C p( |z| ) $,  that the valence of  $ f(z) $  satisfies  
$ N = N(C) < K C \log C $  for some absolute constant  $ K $.
}

We now turn to the three theorems in Section 4.  For the function  $ p(x) $   appearing 
in Theorem 1 we have  $ \tau = 4a(1-a) $,  which means 
by  (\ref{eq:a<1})  that  $ \tau < 1 $  unless  $ a=1/2 $.

 First, consider the case  $ a \neq 1/2 $.   Then  $ \tau < 1 $.  We can, again, 
multiply  $ p(x) $  by a positive  constant  $ C $.  Then for  $ P(x) = Cp(x) $,  
$ \tau $  will be smaller, greater, or equal to 1, and the functions  $ f(z) $  satisfying  
$ Sf(z) = P(z) $   will be of finite or infinite valence in accordance with the 
above discussion.

For  $ a = 1/2 $  we have  $ \tau = 1 $.  It follows from  Theorem 1 that the
function  $ f(z) $  is univalent in the case  $ a = 1/2 $  as well.  But since  $ \tau = 1 $ 
it follows, by  the uniqueness argument discussed earlier,  that the only possible 
case is   $ p(x) = (1 - x^2)^{-2} $  and   $ u(x) = (1 - x^2 )^{1/2} $.   
Theorem 1 supplies exactly the same information.  Indeed,  $ a = 1/2 $ 
implies  $ \lambda = 0 $. (See Figure 2).   Hence we substitute these values in 
(\ref{eq:p-example1})  and get 
\begin{equation*}
	p(x) = \frac{2a + 2a x^2 -4a^2 x^2 }{ (1 - x^2)^2 }  
	     = \frac{1 }{ (1 - x^2)^2 } ,
\end{equation*} 
 as expected.

For the function  $p$ appearing in Theorem 3 the situation is almost identical and 
details are omitted.  On the other hand for Theorem 2 the situation needs more 
care.

By  (\ref{eq:phi(1)-example2}),		
	$ \tau = 4( \lambda + \mu ) (1 -\lambda - \mu) \ge 0 $,
Hence  $ 0 \le \tau < 1 $  unless  $ \lambda + \mu = 1/2 $.  
If  $ \lambda + \mu  \neq  1/2 $,   then $ \tau < 1 $  and the analysis 
is similar to the above.

Now consider the case  $ \lambda + \mu = 1/2 $.  We use the uniqueness argument.
Hence, once again we deduce that there is only one possibility,
namely,  $ u = (1 - x^2)^\lambda \cos^\mu (\pi x/2 ) = (1 - x^2)^{ 1/2 } $.  
This means that necessarily  $ \lambda = 1/2, \ \mu = 0 $.
In fact, this is a somewhat surprising conclusion as it is not included in 
the statement of Theorem 2.   We now explain in more details this point.

Recall that in Theorem 2  the parameter  $ \mu $   may be also negative.  It is easy to 
see that   $ \mu > 0 $   contradicts  $ \lambda + \mu = 1/2 $  and thus can not occur 
with our assumption.   Indeed, this follows at once from
$ (1 - \mu)/2 \le \lambda $,  $ (1 + \mu)/2 \le \lambda + \mu $. 
See (\ref{eq:example2-condition}).      
On the other hand a priori  it is possible to have both  $ \mu < 0 $  and 
$ \lambda + \mu = 1/2 $   without contradicting  (\ref{eq:example2-condition})  
Nevertheless the above analysis shows that the  ``uniqueness argument'' implies 
$ u = (1 - x^2)^\lambda \cos^\mu (\pi x/2 ) = (1 - x^2)^{ 1/2 } $,  i.e., $ \mu = 0 $.

Summing up, we conclude from the above that one may add to the statements of 
Theorem 2 the following additional information:

Suppose  $ \tau =1 $.  Then  $ \mu < 0 $  can not occur, i.e.,  the function arising 
for the values  $ \mu < 0 $  and  $ \lambda + \mu = 1/2 $  is not univalent.  
This means that  $ \varphi(x) = p(x)(1-x^2)^2 $  is not monotone in the interval (0,1). 
Indeed, monotonicity for  $ \varphi $  would supply a function  $ u(x) $ 
$ u = (1 - x^2)^{ 1/2 + |\mu| }  \cos^{ -| \mu | }  (\pi x/2 ) \neq (1 - x^2)^{ 1/2 } $, 
contradicting the uniqueness argument.



\setlength{\parindent}{0 mm}

Department of Mathematics, Technion --- I.I.T., Haifa 32000, Israel	\\ 
{\tt dova@tx.technion.ac.il}	\\ 
{\tt elias@tx.technion.ac.il}

\end{document}